\documentclass[12pt]{article}
\usepackage{amsbsy,amsfonts,amsmath,amssymb,amsthm,bbm,enumerate}
\usepackage{graphicx}

 \numberwithin{equation}{section}
 \theoremstyle{plain}

\newtheorem{theorem}{Theorem}[section]
\newtheorem{corollary}[theorem]{Corollary}

\newtheorem{lemma}[theorem]{Lemma}
\newtheorem{proposition}[theorem]{Proposition}

\newcommand{\bb}{\underline{b}}
\newcommand{\cA}{{\cal A}}
\newcommand{\cC}{{\cal C}}
\newcommand{\cN}{{\cal N}}

\newcommand{\dc}{d_{\rm c}}
\newcommand{\db}{\rightrightarrows}

\newcommand{\dpst}{\displaystyle}

\newcommand{\indic}{{\mathbbm1}}
\newcommand{\ind}[1]{\indic_{\{#1\}}}
\newcommand{\lbeq}[1]{\label{eq:#1}}
\newcommand{\mC}{{\mathbb C}}
\newcommand{\mE}{{\mathbb E}}
\newcommand{\mN}{{\mathbb N}}
\newcommand{\mP}{{\mathbb P}}
\newcommand{\mR}{{\mathbb R}}
\newcommand{\mZ}{{\mathbb Z}}
\newcommand{\nn}{\nonumber}

\newcommand{\pc}{p_{\rm c}}

\newcommand{\Proof}[1]{\paragraph{\it #1.}}
\newcommand{\QED}{\hspace*{\fill}\rule{7pt}{7pt}\smallskip}
\newcommand{\Rd}{\mR^d}

\newcommand{\refeq}[1]{(\ref{eq:#1})}
\newcommand{\sss}{\scriptscriptstyle}

\newcommand{\tb}{\overline{b}}

\newcommand{\Zd}{\mZ^d}
\newcommand{\Zp}{\mZ_+}

\title{Critical behavior and the limit distribution\\
for long-range oriented percolation. I}

\author{
Lung-Chi~Chen\footnote{Department of Mathematics, Fu-Jen Catholic
University, Taiwan.  {\tt lcchen@math.fju.edu.tw}}\\
Akira~Sakai\footnote{Department of Mathematical Sciences, University
of Bath, UK.  {\tt a.sakai@bath.ac.uk}}
\date{March 15, 2007\footnote{Updated: August 14, 2007}}
}

\begin{document}
\maketitle

\begin{abstract}
We consider oriented percolation on $\Zd\times\Zp$ whose bond-occupation
probability is $pD(\,\cdot\,)$, where $p$ is the percolation parameter and $D$
is a probability distribution on $\Zd$.  Suppose that $D(x)$ decays as
$|x|^{-d-\alpha}$ for some $\alpha>0$.  We prove that the two-point function
obeys an infrared bound which implies that various critical exponents take on
their respective mean-field values above the upper-critical dimension
$\dc=2(\alpha\wedge2)$.  We also show that, for every $k$, the Fourier
transform of the normalized two-point function at time $n$, with a proper
spatial scaling, has a convergent subsequence to $e^{-c|k|^{\alpha\wedge2}}$
for some $c>0$.
\end{abstract}

\section{Introduction}
Oriented percolation is a model that exhibits a phase transition when the
percolation parameter $p$ in the bond-occupation probability $pD(\,\cdot\,)$
changes its value, where $D$ is a given probability distribution on $\Zd$.  It
has been proved using the lace expansion \cite{hs02,NY2} that finite-variance
oriented percolation, where the tail of $D$ decays fast enough to ensure finite
variance $\sigma^2=\sum_x|x|^2D(x)$ in particular, exhibits the critical
behavior for (finite-range) branching random walk, if $d>4$ and $\sigma^2\gg1$
or $d\gg4$; it has also been proved that, for every $p\le\pc$ for finite-range
oriented percolation \cite{NY2} and for general (possibly infinite-range)
finite-variance oriented percolation at $p=\pc$ \cite{hs02}, the Fourier
transform of the normalized two-point function at time $n$, spatially scaled by
$\sqrt{n}$, converges to $e^{-c|k|^2}$ for some constant $c\in(0,\infty)$.

In this paper, we consider long-range oriented percolation with index
$\alpha>0$, where $D(x)$ decays as $|x|^{-d-\alpha}$ for large $|x|$. In
\cite{cs06}, Chen and Shieh studied a long-range model with $\alpha=1$ and
proved that, if $d>2$ (and a certain spread-out parameter $L\gg1$), the
standard susceptibility exponent $\gamma$ and a couple of other critical
exponents take on their respective mean-field values.  The goal of this paper
is to investigate the $\alpha$-dependence of the critical behavior and the
limit distribution. We prove that the model exhibits the mean-field behavior if
$d>2(\alpha\wedge2)$ (and a spread-out parameter $L\gg1$).  Furthermore, we
prove that, for every $p\le\pc$, the Fourier transform of the normalized
two-point function at time $n$, spatially scaled by $n^{\frac1{\alpha\wedge2}}$
if $\alpha\ne2$ or by $\sqrt{n\log n}$ if $\alpha=2$, is bounded from below by
$e^{-c|k|^{\alpha\wedge2}}$ and from above by $e^{-c'|k|^{\alpha\wedge2}}$ in
$n\uparrow\infty$, where $c,c'\in(0,\infty)$ and $c/c'=1+O(L^{-d})$.  We stress
that, although we do not prove convergence in this paper, our results hold for
$p\le\pc$ for general finite-variance oriented percolation, which is not
completely covered in the aforementioned results in \cite{hs02,NY2}.

Our proof is based on the lace expansion for oriented percolation.  We analyze
the lace expansion for all $\alpha>0$ simultaneously to discover a potential
crossover in the critical behavior by changing the value of $\alpha$. However,
since our $D$ does not have finite variance when $\alpha\le2$, the standard
Taylor-expansion analyses for the Fourier transform of the expansion
coefficients for finite-variance oriented percolation do not always work.  To
overcome this difficulty, we use the trigonometric techniques that were first
developed in \cite{bchss05} for percolation on finite graphs and later in
\cite{Slad06} for finite-range self-avoiding walk on $\Zd$.  We adapt these
techniques for the time-oriented setting (to analyze the Fourier-Laplace
transform of the expansion coefficients).

\subsection{Model}\label{ss:models}
We define the model more precisely.  A bond is an ordered pair
$((x,n),(y,n+1))$ of vertices in space-time $\Zd\times\Zp$, where
$\Zp\equiv\{0\}\;\Dot\cup\;\mN$ is the set of nonnegative integers.  Each bond
is, independently of the other bonds, occupied (resp., vacant) with probability
$pD(y-x)$ (resp., $1-pD(y-x)$), where $D$ is a probability distribution on
$\Zd$.  The percolation parameter $p\in[0,\|D\|_\infty^{-1}]$ equals the
average number of occupied bonds per vertex.  We say that $(x,l)$ is connected
to $(y,n)$, and write $(x,l)\to(y,n)$, if either $(x,l)=(y,n)$ or there is a
time-oriented path of occupied bonds from $(x,l)$ to $(y,n)$.  Let $\mP_p$ be
the probability distribution of the bond variables, and denote its expectation
by $\mE_p$.

Our $D$ is defined as follows.  Let $h$ be a bounded probability distribution
on $\Rd$ that is invariant under rotations by $\pi/2$ and reflections in the
coordinate hyperplanes.  Suppose that $h$ is piecewise continuous, so that
$\int_{\Rd}d^dx\,h(x)\equiv1$ can be approximated by the Riemann sum
$\frac1{L^d}\sum_{x\in\Zd}h(x/L)$ for large $L<\infty$.  We define
\begin{align}\lbeq{Ddef}
D(x)=\frac{h(x/L)}{\sum_{y\in\Zd}h(y/L)},
\end{align}
where $x/L=(x_1/L,\dots,x_d/L)$.  Note that the denominator is $O(L^d)$.

Fix $\alpha>0$ throughout this paper.  We assume that there is an $\ell<\infty$
such that
\begin{align}\lbeq{majorass}
h(x)\asymp|x|^{-d-\alpha}\qquad(|x|\ge\ell),
\end{align}
where $f(x)\asymp g(x)$ means that $f(x)/g(x)$ is bounded away from zero and
infinity.  We note that the $r^\text{th}$~moment $\sum_{x\in\Zd}|x|^rD(x)$ does
not exist if $r\ge\alpha$, but exists and equals $O(L^r)$ if $r\in(0,\alpha)$.
A simple example of $h$ that satisfies the above assumptions is
\begin{align}\lbeq{heg}
h(x)=\frac1{\cN}\,(|x|\vee1)^{-d-\alpha},
\end{align}
where $\cN$ is the normalization constant.  In this case, $D$ equals
\begin{align}\lbeq{Deg}
D(x)=\frac{(|\frac{x}L|\vee1)^{-d-\alpha}}{\sum_{y\in\Zd}(|\frac{y}L|
 \vee1)^{-d-\alpha}}.
\end{align}

The main properties of $D$ are summarized as follows:

\begin{proposition}\label{prop:Dnsup-bd}
Let $\lambda=L^{-d}$, and denote by $D^{\star n}$ and $\hat D$, respectively,
the $n$-fold convolution and the Fourier transform of $D$:
\begin{align}
D^{\star n}(x)=\begin{cases}
 D(x)&(n=1),\\
 \sum_{y\in\Zd}D^{\star(n-1)}(y)\,D(x-y)&(n\ge2),
 \end{cases}&&&&
\hat D(k)=\sum_{x\in\Zd}e^{ik\cdot x}D(x).
\end{align}
Then, for $L\gg1$, there are $C<\infty$ and $\Delta\in(0,1)$ such that
\begin{align}\lbeq{Dprop}
\|D^{\star n}\|_\infty\le C\lambda\,n^{-\frac{d}{\alpha\wedge2}},&&
 1-\hat D(k)\begin{cases}
 <2-\Delta&(k\in[-\pi,\pi]^d),\\
 >\Delta&(\|k\|_\infty>(\ell L)^{-1}).
 \end{cases}
\end{align}
Moreover, when $\|k\|_\infty\le(\ell L)^{-1}$,
\begin{align}\lbeq{Dprop-k0}
1-\hat D(k)\asymp\begin{cases}
 (L|k|)^{\alpha\wedge2}&(\alpha\ne2),\\
 (L|k|)^2\log\frac\pi{2\ell L|k|}&(\alpha=2).
 \end{cases}
\end{align}
\end{proposition}

We will prove Proposition~\ref{prop:Dnsup-bd} in Appendix~\ref{a:Dprop}.

\subsection{Main results}\label{ss:results}
We investigate the following two-point function:
\begin{align}\lbeq{2pt}
\varphi_p(y-x,n-l)=\mP_p((x,l)\to(y,n)),
\end{align}
where we have used the fact that the right-hand side depends only on $y-x$
and $n-l$, due to the translation invariance of the model.  Assuming
summability of the two-point function, we define, for
$k\in[-\pi,\pi]^d$ and $z\in\mC$,
\begin{align}\lbeq{2pt-FL}
Z_p(k;n)=\sum_{x\in\Zd}\varphi_p(x,n)\,e^{ik\cdot x},&&
\hat\varphi_p(k,z)=\sum_{n\in\Zp}Z_p(k;n)\,z^n.
\end{align}
Let $\cC_n$ be the set of vertices at time $n$ that are connected from $(o,0)$,
and let $\cC=\bigcup_{n\ge0}\cC_n$.  The quantities in \refeq{2pt-FL} for $k=0$
and $(k,z)=(0,1)$ can be described as
\begin{align}
Z_p(0;n)=\mE_p[|\cC_n|],&& \chi_p\equiv\hat\varphi_p(0,1)=\mE_p[|\cC|],
\end{align}
where $|\cA|$ is the cardinality of a set $\cA$, and $\chi_p$ is called the
susceptibility.  Since $Z_p(0;n)$ is sub-multiplicative, i.e., for $l,n\ge0$,
\begin{align}
Z_p(0;l+n)=\sum_{x\in\Zd}\mP_p\bigg(\bigcup_{y\in\Zd}\Big\{\{(o,0)\to(y,l)\}
 \cap\{(y,l)\to(x,l+n)\}\Big\}\bigg)\le Z_p(0;l)\,Z_p(0;n),
\end{align}
the radius $m_p$ of convergence of the series $\hat\varphi_p(0,z)$ is
well-defined and satisfies (cf., e.g., \cite[Appendix II]{G})
\begin{align}\lbeq{mpdef}
m_p^{-1}=\lim_{n\uparrow\infty}Z_p(0;n)^{1/n}=\inf_{n\ge1}Z_p(0;n)^{1/n}.
\end{align}
This implies that $\hat\varphi_p(0,m)$ for $m\in\mR$ diverges as $m\uparrow
m_p$ for every $p>0$, because
\begin{align}\lbeq{chim-ineq}
\hat\varphi_p(0,m)=\sum_{n\in\Zp}Z_p(0;n)\,m^n\ge\sum_{n\in\Zp}
 \bigg(\frac{m}{m_p}\bigg)^n=\frac{m_p}{m_p-m}.
\end{align}
This also implies that $m_p>1$ if and only if $\chi_p<\infty$.  Since
$\hat\varphi_0(0,m)=1$ for any $m\ge0$, we define $m_0=\infty$. It is known
\cite{AB,AN,bg90,gh01} that there is a unique critical point $\pc\ge1$ such
that
\begin{align}
\chi_p\begin{cases}
 <\infty,&\text{if }p<\pc,\\ =\infty,&\text{if }p\ge\pc,
 \end{cases}&&
\Theta_p\equiv\mP_p(|\cC|=\infty)\begin{cases}
 =0,&\text{if }p\le\pc,\\ >0,&\text{if }p>\pc,
 \end{cases}
\end{align}
and that $\lim_{p\uparrow\pc}\chi_p=\infty$ (hence $m_{\pc}\le1$) and
$\lim_{p\downarrow\pc}\Theta_p=0$.

Our first result is about an upper bound on $|\hat\varphi_p(k,z)|$ for $p<\pc$
and $|z|<m_p$.

\begin{theorem}\label{thm:IRbound}
Let $d>2(\alpha\wedge2)$ and $L\gg1$.  Then, there is a $C<\infty$ such that
\begin{align}\lbeq{IRbound}
|\hat\varphi_p(k,z)|\le\frac{C}{p(m_p-|z|)+|\arg(z)|+1-\hat D(k)},
\end{align}
for any $p\in(0,\pc)$, $k\in[-\pi,\pi]^d$ and $z\in\mC$ with $|z|<m_p$.
\end{theorem}

To prove this theorem and the other results throughout this paper, we use the
lace expansion for oriented percolation.  We will briefly review it in
Section~\ref{s:laceexp}.

It has been proved \cite{NY1,NY2} that \refeq{IRbound} holds for
finite-variance oriented percolation (for which, $1-\hat D(k)\asymp|k|^2$) if
$d>4$ and $\sigma^2\gg1$ or $d\gg4$, hence
\begin{align}\lbeq{IRbound-int}
\int_{[-\pi,\pi]^{d+1}}\frac{d^dk}{(2\pi)^d}\frac{d\theta}{2\pi}\,\big|\hat
 \varphi_p(k,me^{i\theta})\big|^3
\end{align}
is bounded uniformly in $p<\pc$ and $m<m_p$.  By the dimension-independent
results in \cite{AN,BA}, this implies that the critical exponents
$\beta,\,\gamma$ and $\delta$ defined as
\begin{align}
\Theta_p\underset{p\downarrow\pc}\asymp(p-\pc)^\beta,&&
 \chi_p\underset{p\uparrow\pc}\asymp(\pc-p)^{-\gamma},&&
 \mP_{\pc}(|\cC|\ge n)\underset{n\uparrow\infty}\asymp n^{-1/\delta},
\end{align}
exist and take on their mean-field values for $d>4$: $\beta=\gamma=1$ and
$\delta=2$. Since our $1-\hat D(k)$ satisfies \refeq{Dprop-k0}, the integral
\refeq{IRbound-int} is bounded uniformly in $p<\pc$ and $m<m_p$ when
$d>2(\alpha\wedge2)$.  Let $\tau$ and $\eta$ be the critical exponents for
$m_p-m_{\pc}$ and $Z_{\pc}(0;n)$, respectively:
\begin{align}
m_p-m_{\pc}\underset{p\uparrow\pc}\asymp(\pc-p)^\tau,&&
Z_{\pc}(0;n)\underset{n\uparrow\infty}\asymp n^\eta.
\end{align}

\begin{corollary}\label{cor:bgdz}
Let $d>2(\alpha\wedge2)$ and $L\gg1$, so that Theorem~\ref{thm:IRbound} holds.
Then, $m_{\pc}=1$ and the critical exponents $\beta,\gamma,\delta$ and $\tau$
exist and take on their respective mean-field values: $\beta=\gamma=\tau=1$ and
$\delta=2$.
\end{corollary}

The identity $\tau=1$ follows immediately from $\gamma=1$ and the inequality
\begin{align}\lbeq{mvschi}
\frac{m_p}{\chi_p}\le m_p-1\le\frac{C}{p\chi_p}\qquad(0<p<\pc).
\end{align}
The lower bound is due to \refeq{chim-ineq} for $m=1$, and the upper bound is
due to Theorem~\ref{thm:IRbound} for $(k,z)=(0,1)$.  By the continuity of
$\chi_p^{-1}$ in $p$, we obtain $m_{\pc}=\lim_{p\uparrow\pc}m_p=1$.  It may be
worth pointing out that the trivial bound $Z_p(0;n)\le p^n$ and the inequality
\refeq{mvschi} with $\chi_p\ge1$ imply $m_p\asymp p^{-1}$ for all $p\in(0,1)$.

The mean-field result on the exponent $\eta$ is in Theorem~\ref{thm:limiting}
below.

The critical exponents are generally believed to be universal in the sense that
their values depend only on $d$ and $\alpha$, but not on the microscopic
details of the model, such as the value of $L<\infty$.  However, the value of
$\pc$ is not universal and changes depending on the value of $L$.
In~\cite{vdHS05}, an asymptotic estimate of $\pc$ as $L\to\infty$ was
investigated for various finite-variance models, such as self-avoiding walk,
percolation, oriented percolation and the contact process, above the
model-dependent upper-critical dimension.  Using
Proposition~\ref{prop:Dnsup-bd} and Theorem~\ref{thm:IRbound}, we obtain the
same asymptotic estimate of $\pc$ for our long-range oriented percolation for
$d>2(\alpha\wedge2)$, as follows:

\begin{theorem}\label{thm:pc}
Let $d>2(\alpha\wedge2)$.  Then, as $L\to\infty$,
\begin{align}\lbeq{pc}
\pc=1+\frac12\sum_{n=2}^\infty D^{\star2n}(o)+O(\lambda^2),
\end{align}
where the sum of the $2n$-fold convolutions over $n\ge2$ is $O(\lambda)$ if
$d>\alpha\wedge2$.
\end{theorem}

Our last results are about asymptotic estimates of the expected number
$Z_p(0;n)$ of vertices at time $n$ connected from $(o,0)$ and the Fourier
transform of the normalized two-point function $Z_p(\,\cdot\,;n)/Z_p(0;n)$. For
finite-range oriented percolation with $d>4$ and $\sigma^2\gg1$ or $d\gg4$,
Nguyen and Yang \cite{NY2} used Tauberian estimates to prove that, for any
$p\in(0,\pc]$ and $k\in\mR^d$, there are $c_1,c_2=1+O(\lambda)$ such that
$Z_p(0;n)\sim c_1m_p^{-n}$ and $Z_p(k/\sqrt{n};n)/Z_p(0;n)\sim e^{-c_2|k|^2}$;
sharper error estimates for general finite-variance oriented percolation at
$p=\pc$ were obtained in \cite{hs02} by an inductive analysis of the lace
expansion.  In this paper, we follow the line of \cite{NY2} using Tauberian
estimates to prove the following theorem for long-range oriented percolation:

\begin{theorem}\label{thm:limiting}
Let $d>2(\alpha\wedge2)$ and $L\gg1$, so that Theorem~\ref{thm:IRbound} holds.
Fix $\epsilon\in(0,1\wedge\frac{d-2(\alpha\wedge2)}{\alpha\wedge2})$.  Then,
the following (i)--(ii) hold for any $p\in(0,\pc]$ and $k\in\Rd$:
\begin{enumerate}[(i)]
\item
There is a $C_1=1+O(\lambda)$ such that
\begin{align}
Z_p(0;n)=C_1m_p^{-n}\big(1+O(n^{-\epsilon})\big)\qquad(n\ge1).
\end{align}
In particular, the critical exponent $\eta$ takes on its mean-field
value: $\eta=0$.
\item
Suppose that there is an $L$-dependent constant $v_\alpha\in(0,\infty)$ such
that
\begin{align}\lbeq{valpha}
1-\hat D(k)\underset{|k|\to0}\sim\begin{cases}
 v_\alpha|k|^{\alpha\wedge2}&(\alpha\ne2),\\
 v_2|k|^2\log\frac1{|k|}&(\alpha=2).
\end{cases}
\end{align}
Let
\begin{align}\lbeq{kn-def}
k_n=k\times
\begin{cases}
(v_\alpha n)^{-\frac1{\alpha\wedge2}}&(\alpha\ne2),\\
(v_2n\log\sqrt{n})^{-\frac12}&(\alpha=2).
\end{cases}
\end{align}
Then, there are $C_2$ and $C_2'$, both equal to $1+O(\lambda)$, such that
\begin{align}\lbeq{Zlim}
e^{-C_2|k|^{\alpha\wedge2}}\le\liminf_{n\to\infty}\frac{Z_p(k_n;n)}{Z_p
 (0;n)}\le\limsup_{n\to\infty}\frac{Z_p(k_n;n)}{Z_p(0;n)}\le e^{-C_2'
 |k|^{\alpha\wedge2}}.
\end{align}
\end{enumerate}
\end{theorem}

We note that our $D$ satisfies the bound \refeq{Dprop-k0} on $1-\hat D(k)$ for
small $k$.  The assumption \refeq{valpha} identifies the coefficient of the
leading term of $1-\hat D(k)$.

In the proof of the above theorem, we estimate fractional moments for the
\emph{time} variable of the lace-expansion coefficients.  In the ongoing work
\cite{cs??}, we have been able to show that the limit of $Z_p(k_n;n)/Z_p(0;n)$
exists for $\alpha>2$ and $d>6$ by crude fractional-moment estimates for the
\emph{spatial} variable of the expansion coefficients.  The difficulty in
proving existence of the limit for all $\alpha>0$ and $d>2(\alpha\wedge2)$ is
due to the fact that the support of our $D$ is unbounded, so that we cannot
simply bound $|x|^r\varphi_p(x,n)$ for some $r>0$, which may show up in the
fractional-moment analysis, by a multiple of $n^r\varphi_p(x,n)$, as done in
\cite{NY2} for finite-range oriented percolation.  To squeeze the bounds in
\refeq{Zlim} in order to identify the limit of $Z_p(k_n;n)/Z_p(0;n)$, we may
have to improve the aforementioned fractional-moment estimates for the
\emph{spatial} variable.  We expect that the idea may also be extended to
investigate $\xi_p^{\sss(r)}(n)\equiv\sum_x|x|^r\varphi_p(x,n)/Z_p(0;n)$.
Nguyen and Yang proved in \cite{NY2} that $\xi_p^{\sss(2)}(n)\asymp n$ for any
$p\in(0,\pc]$ for sufficiently spread-out finite-range oriented percolation for
$d>4$.  We are aiming to show that $\xi_p^{\sss(r)}(n)\asymp
n^{\sss\frac{r}{\alpha\wedge2}}$ for any $p\in(0,\pc]$ and $r<\alpha$ for our
long-range oriented percolation for $d>2(\alpha\wedge2)$.

\subsection{Organization}
The rest of this paper is organized as follows.  In
Section~\ref{s:proofoftheorems}, we prove the above three theorems assuming a
couple of key propositions.  These propositions are proved in
Sections~\ref{s:bootstrapping}--\ref{s:tauber}.  Finally, in the Appendix, we
prove Proposition~\ref{prop:Dnsup-bd}.

\section{Proof of the main results}\label{s:proofoftheorems}
In Sections~\ref{ss:IRbound}--\ref{ss:limiting}, we prove
Theorems~\ref{thm:IRbound}, \ref{thm:pc} and \ref{thm:limiting}, respectively,
assuming several key ingredients.  The most important ingredient is the lace
expansion.

\subsection{Lace expansion}\label{ss:lace-exp}
The idea of the lace expansion was initiated by Brydges and Spencer in
\cite{BS} for investigating weakly self-avoiding walk for $d>4$.  Later, the
lace expansion was applied to various stochastic-geometrical models, such as
strictly self-avoiding walk for $d>4$ (e.g., \cite{HS1}), lattice trees/animals
for $d>8$ (e.g., \cite{HS4}), percolation for $d>6$ (e.g., \cite{HS3}),
oriented percolation for $d>4$ (e.g., \cite{NY1}) and the contact process for
$d>4$ (e.g., \cite{S1}).  Application to the Ising model was recently reported
in \cite{s06}.  See \cite{Slad06} for a complete list of references up to 2005.

The derivation of the lace expansion, the definition of the expansion
coefficients and their diagrammatic bounds in terms of two-point functions
depend on which model is concerned, but are independent of the specific choice
of $D$.  Therefore, we can apply the standard lace expansion for oriented
percolation to the current long-range setting.  We will briefly review the
expansion in Section~\ref{s:laceexp}.

The result of the lace expansion is a recursion equation similar
to that for the random-walk two-point function
\begin{align}\lbeq{laceexp-rw}
P_p(x,n)=\delta_{x,o}\delta_{n,0}+p^nD^{\star n}(x)\ind{n\ge1}
 =\delta_{x,o}\delta_{n,0}+(q_p*P_p)(x,n),
\end{align}
where $\ind{\cdots}$ is the indicator function and
\begin{align}\lbeq{q-def}
q_p(x,n)=pD(x)\delta_{n,1}.
\end{align}
For oriented percolation, we have (see Proposition~\ref{prop:laceexp-2pt}
below)
\begin{align}\lbeq{laceexp-x}
\varphi_p(x,n)=\pi_p(x,n)+(\pi_p*q_p*\varphi_p)(x,n)\qquad
 (0\le p\le\|D\|_\infty^{-1}),
\end{align}
where $\pi_p(x,n)$ is the alternating sum of the nonnegative lace-expansion
coefficients $\pi_p^{\sss(N)}(x,n)$:
\begin{align}
\pi_p^{\sss(N)}(x,n)\ge0\qquad(N=0,1,\dots),&&
\pi_p(x,n)=\sum_{N=0}^\infty(-1)^N\pi_p^{\sss(N)}(x,n).
\end{align}
If $n=0$, then $\pi_p^{\sss(N)}(x,0)=\delta_{x,o}\delta_{N,0}$, hence
$\pi_p(x,0)=\delta_{x,o}$, due to the definition \refeq{piN-def} of
$\pi_p^{\sss(N)}(x,n)$ below.  Comparing \refeq{laceexp-rw} and
\refeq{laceexp-x}, we are naturally led to expect that $\varphi_p(x,n)$ behaves
similarly to $P_p(x,n)$, if $\pi_p(x,n)-\delta_{x,o}\delta_{n,0}$ is small.

\subsection{Infrared bound}\label{ss:IRbound}
We prove Theorem~\ref{thm:IRbound} by comparing $\hat\varphi_p(k,z)$,
where $k\in[-\pi,\pi]^d$ and $z\in\mC$ with $|z|<m_p$, with the Fourier
transform of the random-walk Green's function with a certain rate
$\mu=\mu_p(z)\in\mC$:
\begin{align}\lbeq{G-Fourier}
\hat G_\mu(k)\equiv\sum_{(x,n)\in\Zd\times\Zp}P_\mu(x,n)\,e^{ik\cdot x}
 =\frac1{1-\mu\hat D(k)}\qquad(|\mu\hat D(k)|<1).
\end{align}
It is not hard to see that $\hat G_\mu(k)$ obeys the following infrared bound:
\begin{align}\lbeq{G-IRbd}
|\hat G_\mu(k)|\le\frac{c}{(1-|\mu|)+|\arg(\mu)|+1-\hat D(k)},
\end{align}
where $c<\infty$ is independent of $\mu$ and $k$.

Let
\begin{align}\lbeq{mu-def}
\mu_p(z)=\big(1-\hat\varphi_p(0,|z|)^{-1}\big)\,e^{i\arg(z)},
\end{align}
where $|\mu_p(z)|<1$ for $|z|<m_p$ and $\mu_p(m)\uparrow1$ as $m\uparrow m_p$.
Inspired by the bootstrapping hypotheses used in \cite{bchss05} for percolation
on finite graphs and in \cite{Slad06} for finite-range self-avoiding walk on
$\Zd$, we define
\begin{align}
f(p,m)=\max_{i=1,2,3}f_i(p,m)\qquad(p<\pc,~m<m_p),
\end{align}
where
\begin{gather}
 f_1(p,m)=p(m\vee1),\qquad\qquad
 f_2(p,m)=\sup_{\substack{k\in[-\pi,\pi]^d\\ z\in\mC:|z|\in\{m,1\}}}\bigg|
  \frac{\hat\varphi_p(k,z)}{\hat G_{\mu_p(z)}(k)}\bigg|,\qquad\qquad
  \lbeq{f1f2-def}\\
 f_3(p,m)=\sup_{\substack{k,l\in[-\pi,\pi]^d\\ z\in\mC:|z|\in\{m,1\}}}\frac{
  \hat G_{\mu_p(m\vee1)}(k)~|\hat\varphi_p(l,z)-\frac12(\hat\varphi_p(l+k,z)
  +\hat\varphi_p(l-k,z))|}{K\sum_{(j,j')=(0,\pm1),(1,-1)}|\hat G_{\mu_p(z)}
  (l+jk)\,\hat G_{\mu_p(z)}(l+j'k)|}.\lbeq{f3-def}
\end{gather}
for some large but finite constant $K>0$ whose precise value is unimportant for
the moment and will be determined in Section~\ref{ss:IR(ii)}.  These functions
will be used in the bootstrapping argument, as stated in
Proposition~\ref{prop:bootstrapping} below.  We emphasize that, although the
work in \cite{bchss05,Slad06} did not concern the long-range models, the
definition of $f_3$ is well-adapted to the long-range setting, especially for
$\alpha\le2$; since we are not using the Taylor expansion for the numerator of
\refeq{f3-def}, we do not have to assume convergence of the second moment for
the spatial variable of the two-point function.  We use similar functions in
the bootstrapping argument in \cite{hhs?} to investigate the critical behavior
for the long-range Ising model, percolation and self-avoiding walk on $\Zd$.

We prove below Theorem~\ref{thm:IRbound} using the following proposition:

\begin{proposition}\label{prop:bootstrapping}
\begin{enumerate}[(i)]
\item
Let $d>2(\alpha\wedge2)$ and $L\gg1$ and fix $p<\pc$ and $m<m_p$.  Then,
$f(p,m)\le3$ implies that there is a $(p,m)$-independent constant $C<\infty$
such that
\begin{align}
\sum_{(x,n)\in\Zd\times\mN}n^r\pi_p^{\sss(N)}(x,n)m^n&\le(C\lambda)^{N\vee1}
 \qquad(N\ge0,~r=0,1),\lbeq{piN-bd}\\
\sum_{(x,n)\in\Zd\times\Zp}\big(1-\cos(k\cdot x)\big)|\pi_p(x,n)|m^n
 &\le C\lambda\,\hat G_{\mu_p(m\vee1)}(k)^{-1}\qquad(k\in[-\pi,\pi]^d).
 \lbeq{piNcos-bd}
\end{align}
\item
Let $d>2(\alpha\wedge2)$ and $L\gg1$ and fix $p<\pc$ and $m<m_p$.  Then,
\refeq{piN-bd}--\refeq{piNcos-bd} and $f(p,m)\le3$ imply the stronger bound
$f(p,m)\le2$.
\item
The function $f(p,m)$ is continuous in $m<m_p$ for every $p<\pc$, and $f(p,1)$
is continuous in $p<\pc$, with $f(0,1)=1$.
\end{enumerate}
\end{proposition}

We will prove Proposition~\ref{prop:bootstrapping} in
Section~\ref{s:bootstrapping}.

\Proof{Proof of Theorem~\ref{thm:IRbound} assuming
 Proposition~\ref{prop:bootstrapping}}
Note that Proposition~\ref{prop:bootstrapping}(i)--(ii) imply
$f(p,m)\notin[2,3)$ for every $p<\pc$ and $m<m_p$.  With the help of the
continuity in Proposition~\ref{prop:bootstrapping}(iii), we conclude that
indeed $f(p,m)\le2$ holds for all $p<\pc$ and $m<m_p$.  In particular, by
\refeq{G-IRbd} and the definition of $f_2$, we have
\begin{align}\lbeq{IRprebd}
|\hat\varphi_p(k,z)|\le\frac{2c}{(1-|\mu_p(z)|)+|\arg(z)|+1-\hat D(k)}
 \qquad(p<\pc,~|z|<m_p).
\end{align}
To complete the proof of Theorem~\ref{thm:IRbound}, it suffices to show that
\begin{align}\lbeq{1-mu-bd}
1-|\mu_p(z)|\equiv\hat\varphi_p(0,|z|)^{-1}\ge\frac12p(m_p-|z|)
 \qquad(0<p<\pc).
\end{align}

Before proving \refeq{1-mu-bd}, we note that $\hat\varphi_p(0,m)$ diverges as
$m\uparrow m_p$ for every $p>0$ (cf., \refeq{chim-ineq}) and that, by using
\refeq{laceexp-x},
\begin{align}\lbeq{chim}
1\le\sum_{(x,n)\in\Zd\times\Zp}\varphi_p(x,n)\,m^n=\hat\varphi_p(0,m)
 =\frac{\hat\pi_p(0,m)}{1-pm\hat\pi_p(0,m)}<\infty\qquad(m<m_p).
\end{align}
By \refeq{piN-bd} for $r=0$, $|\hat\pi_p(0,m)-1|$ is uniformly bounded by
$O(\lambda)$.  Moreover, by monotone convergence and \refeq{piN-bd} for $r=1$,
\begin{align}\lbeq{MCT2}
m_p|\hat\pi_p(0,m_p)-\hat\pi_p(0,m)|&\le\sum_{(x,n)}|\pi_p(x,n)|\,m_p(m_p^n
 -m^n)\nn\\
&\le(m_p-m)\sum_{(x,n)}n|\pi_p(x,n)|\,m_p^n\nn\\
&\le(m_p-m)\sum_{(x,n)}\sum_{N=0}^\infty n\,\pi_p^{\sss(N)}(x,n)\,m_p^n\nn\\
&=(m_p-m)\lim_{m\uparrow m_p}\sum_{(x,n)}\sum_{N=0}^\infty n\,\pi_p^{\sss
 (N)}(x,n)\,m^n\le O(\lambda)(m_p-m),
\end{align}
where the $O(\lambda)$ term is independent of $m$, so that
$\hat\pi_p(0,m_p)=\lim_{m\uparrow m_p}\hat\pi_p(0,m)$.  Therefore, for
$\hat\varphi_p(0,m)$ to diverge as $m\uparrow m_p$, the denominator in
\refeq{chim} should be nonnegative and vanish as $m\uparrow m_p$, and hence
\begin{align}\lbeq{pmppi}
pm_p\hat\pi_p(0,m_p)=1\qquad(0<p<\pc).
\end{align}

Now we continue with the proof of \refeq{1-mu-bd}.  Since
$\hat\pi_p(0,|z|)=1+O(\lambda)>0$ as explained above, we obtain
\begin{align}
\hat\varphi_p(0,|z|)^{-1}=\bigg(\frac{\hat\pi_p(0,|z|)}{1-p|z|
 \hat\pi_p(0,|z|)}\bigg)^{-1}=\hat\pi_p(0,|z|)^{-1}-p|z|.
\end{align}
By repeated use of \refeq{pmppi}, $\hat\varphi_p(0,|z|)^{-1}$ is
rewritten as
\begin{align}
\hat\varphi_p(0,|z|)^{-1}&=\hat\pi_p(0,|z|)^{-1}-p|z|+pm_p-\hat
 \pi_p(0,m_p)^{-1}=p(m_p-|z|)+\frac{\hat\pi_p(0,m_p)-\hat\pi_p
 (0,|z|)}{\hat\pi_p(0,|z|)\,\hat\pi_p(0,m_p)}\nn\\
&=p\bigg((m_p-|z|)+\frac{m_p(\hat\pi_p(0,m_p)-\hat\pi_p(0,|z|))}
 {\hat\pi_p(0,|z|)}\bigg).
\end{align}
By \refeq{MCT2}, we have arrived at
\begin{align}
\hat\varphi_p(0,|z|)^{-1}\ge\big(1-O(\lambda)\big)p(m_p-|z|).
\end{align}
This completes the proof of Theorem~\ref{thm:IRbound} assuming
Proposition~\ref{prop:bootstrapping}. \QED

\subsection{Asymptotic estimate of $\pc$}\label{ss:pc}
We begin with the identity \refeq{chim} for $m=1$:
\begin{align}\lbeq{chim=1}
1\le\chi_p\equiv\hat\varphi_p(0,1)=\frac{\hat\pi_p(0,1)}{1-p\hat\pi_p
 (0,1)}<\infty\qquad(p<\pc).
\end{align}
By \refeq{piN-bd} for $m=1$ and $r=0$, $|\hat\pi_p(0,1)-1|$ is bounded by
$O(\lambda)$ uniformly in $p<\pc$.  Since $\chi_p\uparrow\infty$ and
$m_p\downarrow1$ as $p\uparrow\pc$, we have
\begin{align}\lbeq{pcpi}
1=\pc\hat\pi_{\pc}(0,1)\equiv\pc\lim_{p\uparrow\pc}\hat\pi_p(0,1),
\end{align}
and therefore $\pc=\hat\pi_{\pc}(0,1)^{-1}=1+O(\lambda)$.

To improve this estimate, we use the following proposition:

\begin{proposition}\label{prop:piderbd}
Let $d>2(\alpha\wedge2)$ and $L\gg1$.  Then, there is a $C<\infty$ such that,
for $p\in(1,\pc)$,
\begin{align}\lbeq{piderbd}
|\partial_p\hat\pi_p(0,1)|\le C\lambda.
\end{align}
\end{proposition}

We will prove Proposition~\ref{prop:piderbd} in Section~\ref{s:pc}.

\Proof{Proof of Theorem~\ref{thm:pc} assuming Proposition~\ref{prop:piderbd}}
First we rewrite \refeq{pcpi} as
\begin{align}\lbeq{pcpi-rewr}
1=\pc\big(\hat\pi_{\pc}(0,1)-\hat\pi_1(0,1)\big)+(\pc-1)\,\big(\hat\pi_1(0,1)
 -1\big)+\big(\hat\pi_1(0,1)-1\big)+\pc.
\end{align}
We already know $(\pc-1)(\hat\pi_1(0,1)-1)=O(\lambda^2)$.  By the mean-value
theorem and Proposition~\ref{prop:piderbd},
\begin{align}
|\hat\pi_{\pc}(0,1)-\hat\pi_1(0,1)|=(\pc-1)|\partial_p\hat\pi_p(0,1)|\le
 O(\lambda^2).
\end{align}
Moreover, by \refeq{piN-bd} for $(p,m)=(1,1)$ and $r=0$, we have
$\hat\pi_1^{\sss(N)}(0,1)\le O(\lambda)^N$ for $N\ge2$.  Therefore,
\begin{align}
\pc=1+\hat\pi_1^{\sss(1)}(0,1)-\big(\hat\pi_1^{\sss(0)}(0,1)-1\big)
 +O(\lambda^2).
\end{align}

To complete the proof of Theorem~\ref{thm:pc}, it suffices to show that,
for $d>2(\alpha\wedge2)$,
\begin{align}\lbeq{suff-pcasy}
\hat\pi_1^{\sss(1)}(0,1)-\big(\hat\pi_1^{\sss(0)}(0,1)-1\big)=\frac12
 \sum_{n=2}^\infty D^{\star2n}(o)+O(\lambda^2),
\end{align}
where the sum is $O(\lambda)$ if $d>\alpha\wedge2$, because of
Proposition~\ref{prop:Dnsup-bd}.  In fact, \refeq{suff-pcasy} follows from the
same argument as in \cite[Section~3.1]{vdHS05} and using
Proposition~\ref{prop:Dnsup-bd}.  The main point is that, since $p=1$, we can
estimate $\hat\pi_1^{\sss(i)}(0,1)$ with random walks.  For example,
$\hat\pi_p^{\sss(0)}(0,1)-1$ is the sum over $(x,n)\in\Zd\times\mN$ of the
probability that there are at least \textit{two} bond-disjoint connections from
$(o,0)$ to $(x,n)$ (cf., the definition \refeq{pi0-def} of
$\pi_p^{\sss(0)}(x,n)$ below).  Since $p=1$, each of these bond-disjoint
connections can be approximated by a random-walk path from $o$ to $x$ in $n$
steps.  Therefore, the main contribution to $\hat\pi_1^{\sss(0)}(0,1)-1$ is
\begin{align}
\frac12\sum_{n=2}^\infty\sum_{x\in\Zd}\big(D^{\star n}(x)\big)^2=\frac12
 \sum_{n=2}^\infty D^{\star2n}(o),
\end{align}
where the combinatorial factor $\frac12$ is due to the symmetry between the two
bond-disjoint connections (cf., \cite[(3.11)]{vdHS05}), which is absent in the
main contribution to $\hat\pi_1^{\sss(1)}(0,1)$ (cf., \cite[(3.22)]{vdHS05}),
leading to the factor $\frac12$ in the difference \refeq{suff-pcasy}.  The
corrections to $\hat\pi_1^{\sss(0)}(0,1)-1$ and $\hat\pi_1^{\sss(1)}(0,1)$ can
be estimated as $O(\lambda^2)$ by applying Proposition~\ref{prop:Dnsup-bd} to
the error terms in \cite[Section~3.1]{vdHS05}.  For example,
\cite[(3.29)]{vdHS05} is replaced by
\begin{align}
\sum_{\substack{t,s,s'\in\Zp:\\ 0\le s<s'\le t}}\frac{O(\lambda)}{(1\vee
 t)^{d/(\alpha\wedge2)}}\,\frac{O(\lambda)}{(1\vee(s'-s))^{d/(\alpha\wedge
 2)}}\le\sum_{t=0}^\infty\frac{O(\lambda^2)}{(1\vee t)^{d/(\alpha\wedge2)
 -1}}\le O(\lambda^2),
\end{align}
where we have used $d>2(\alpha\wedge2)$.  This completes the proof of
Theorem~\ref{thm:pc} assuming Proposition~\ref{prop:piderbd}. \QED

\subsection{Limit distribution}\label{ss:limiting}
Assuming the lace expansion \refeq{laceexp-x} and the bounds in
Proposition~\ref{prop:bootstrapping} on the expansion coefficients, we have
that, for $p\in(0,\pc)$, $k\in[-\pi,\pi]^d$ and $m<m_p$,
\begin{align}\lbeq{reorg0}
\hat\varphi_p(k,m)^{-1}=\hat\pi_p(k,m)^{-1}-pm\hat D(k),
\end{align}
where $\hat\pi_p(k,m)=1+O(\lambda)$.  In the course of the proof of
Theorem~\ref{thm:IRbound} in Section~\ref{ss:IRbound}, we obtained
$pm_p\equiv\hat\pi_p(0,m_p)^{-1}=1+O(\lambda)$ for $p\in(0,\pc]$ and
$m_{\pc}=1$, as stated in Corollary~\ref{cor:bgdz}.  For $m<1$,
$\hat\pi_{\pc}(k,m)\equiv\lim_{p\uparrow\pc}\hat\pi_p(k,m)$ is well-defined,
due to \refeq{reorg0} and the continuity of $\hat\varphi_p(k,m)$ in $p<\pc$ for
every $m<1$, as well as the uniform bound on $\hat\pi_p(k,m)$.

Using these facts and Tauberian estimates, we first derive an asymptotic
formula of $Z_p(k;n)$ for every $p\in(0,\pc]$.  Then, by using this formula, we
will prove Theorem~\ref{thm:limiting}.

Since $\hat\pi_p(k,m)=1+O(\lambda)$ and $\hat\pi_p(0,m_p)^{-1}=pm_p$, we can
reorganize \refeq{reorg0} for $m<m_p$ as
\begin{align}\lbeq{reorg1}
\hat\varphi_p(k,m)^{-1}&=\underbrace{\hat\pi_p(k,m)^{-1}-pm\hat D(k)-\Big(
 \hat\pi_p(k,m_p)^{-1}-pm_p\hat D(k)\Big)}_{p(m_p-m)\hat A_p(k,m)}\nn\\
&\quad+\underbrace{\hat\pi_p(k,m_p)^{-1}-pm_p\hat D(k)-\Big(\hat\pi_p(0,
 m_p)^{-1}-pm_p\Big)}_{pm_p\hat B_p(k)}\nn\\
&=pm_p\Big(\big(1-\tfrac{m}{m_p}\big)\hat A_p(k,m)+\hat B_p(k)\Big),
\end{align}
where
\begin{align}
\hat A_p(k,m)&=\hat D(k)-\frac{\hat\pi_p(k,m_p)^{-1}-\hat\pi_p(k,m)^{-1}}
 {p(m_p-m)},\\
\hat B_p(k)&=1-\hat D(k)+\frac{\hat\pi_p(k,m_p)^{-1}-\hat\pi_p
 (0,m_p)^{-1}}{pm_p}.\lbeq{B-def}
\end{align}
Similarly to \refeq{MCT2}, we can show that the second term in $\hat A_p(k,m)$
is $O(\lambda)$ and the last term in $\hat B_p(k)$ is $O(\lambda)\hat
G_{\mu_p(m_p\vee1)}(k)\equiv O(\lambda)(1-\hat D(k))$ for $p\le\pc$,
$k\in[-\pi,\pi]^d$ and $m<m_p$.  Then, we decompose $\hat A_p(k,m)$ as $\hat
A_p(k,m)=\hat A_p^{\sss(1)}(k)+\hat A_p^{\sss(2)}(k,m)$, where
\begin{align}
\hat A_p^{\sss(1)}(k)&=\hat D(k)-\frac{m_p\,\partial_m\hat\pi_p(k,m_p)^{-1}}
 {pm_p},\lbeq{A1-def}\\
\hat A_p^{\sss(2)}(k,m)&=\frac{m_p\,\partial_m\hat\pi_p(k,m_p)^{-1}}{pm_p}
 -\frac{\hat\pi_p(k,m_p)^{-1}-\hat\pi_p(k,m)^{-1}}{p(m_p-m)},\lbeq{A2-def}
\end{align}
where $\partial_m\hat\pi_p(k,m_p)^{-1}$ is an abbreviation for
$\partial_m\hat\pi_p(k,m)^{-1}|_{m=m_p}$.  Again, similarly to \refeq{MCT2}, we
can show that the common term in \refeq{A1-def}--\refeq{A2-def} is $O(\lambda)$
for any $p\le\pc$ and $k\in[-\pi,\pi]^d$.  In particular, $\hat
A_p^{\sss(1)}(k)$ is continuous at $k=0$, and $\hat A_p^{\sss(1)}(k)+\hat
B_p(k)=1+O(\lambda)$.  Using these quantities, we can rewrite \refeq{reorg1} as
\begin{align}\lbeq{reorg2}
pm_p\hat\varphi_p(k,m)=\frac1{(1-\frac{m}{m_p})\hat A_p(k,m)+\hat B_p(k)}=
 \frac1{(1-\frac{m}{m_p})\hat A_p^{\sss(1)}(k)+\hat B_p(k)}+\hat\Phi_p(k,m),
\end{align}
where
\begin{align}\lbeq{psi-def}
\hat\Phi_p(k,m)=\frac{-(1-\tfrac{m}{m_p})\hat A_p^{\sss(2)}(k,m)}{\big(
 (1-\frac{m}{m_p})\hat A_p(k,m)+\hat B_p(k)\big)\big((1-\frac{m}{m_p})\hat
 A_p^{\sss(1)}(k)+\hat B_p(k)\big)}.
\end{align}
The first term of the rightmost expression in \refeq{reorg2} can be
expanded in powers of $\frac{m}{m_p}$ as
\begin{align}\lbeq{series-exp}
\frac1{\hat A_p^{\sss(1)}(k)+\hat B_p(k)-\frac{m}{m_p}\hat A_p^{\sss(1)}
 (k)}=\frac1{\hat A_p^{\sss(1)}(k)+\hat B_p(k)}\sum_{n=0}^\infty\bigg(
 \frac{m}{m_p}\bigg)^n\bigg(\frac{\hat A_p^{\sss(1)}(k)}{\hat A_p^{\sss
 (1)}(k)+\hat B_p(k)}\bigg)^n.
\end{align}
In Section~\ref{s:tauber}, we will prove the following bound on
$\hat\Phi_p(k,m)$:

\begin{proposition}\label{prp:tauber-suff}
Let $d>2(\alpha\wedge2)$ and $L\gg1$, and fix an
$\epsilon\in(0,1\wedge\frac{d-2(\alpha\wedge2)}{\alpha\wedge2})$.  Then, there
is an $\epsilon$-dependent constant $C_\epsilon<\infty$ such that
\begin{align}\lbeq{psi-est}
|\partial_\zeta\hat\Phi_p(k,m_p\zeta)|\le C_\epsilon|1-\zeta|^{-2+\epsilon}
\end{align}
holds for $p\in(0,\pc]$, $k\in[-\pi,\pi]^d$ and $\zeta\in\mC$ with $|\zeta|<1$.
\end{proposition}

By this result and \cite[Lemma~6.3.3(ii)]{ms93}, the coefficient of
$\zeta^n\equiv(\frac{m}{m_p})^n$ in $\hat\Phi_p(k,m)$ is bounded by
$O(n^{-\epsilon'})$ for any $\epsilon'<\epsilon$.  Together with \refeq{reorg2}
and \refeq{series-exp} and using $pm_p=1+O(\lambda)$, we finally obtain
\begin{align}\lbeq{limit-main}
Z_p(k;n)=\frac{m_p^{-n}}{pm_p(\hat A_p^{\sss(1)}(k)+\hat B_p(k))}\bigg(
 \frac{\hat A_p^{\sss(1)}(k)}{\hat A_p^{\sss(1)}(k)+\hat B_p(k)}\bigg)^n
 +O(m_p^{-n}n^{-\epsilon'})
 \qquad(n\ge1).
\end{align}

\Proof{Proof of Theorem~\ref{thm:limiting} using \refeq{limit-main}} When
$k=0$, since $\hat B_p(0)\equiv0$, we immediately obtain from
\refeq{limit-main} that
\begin{align}
Z_p(0;n)=C_1m_p^{-n}+O(m_p^{-n}n^{-\epsilon'})
 \qquad(n\ge1),
\end{align}
where $C_1\equiv(pm_p\hat A_p^{\sss(1)}(0))^{-1}=1+O(\lambda)$.  This
completes the proof of Theorem~\ref{thm:limiting}(i).

To prove Theorem~\ref{thm:limiting}(ii) using \refeq{limit-main}, it
suffices to investigate
\begin{align}\lbeq{Psi-rewr}
\bigg(\frac{\hat A_p^{\sss(1)}(k)}{\hat A_p^{\sss(1)}
 (k)+\hat B_p(k)}\bigg)^n=\Bigg(\bigg(1+\frac{\hat B_p(k)}{\hat A_p^{\sss
 (1)}(k)}\bigg)^{\frac{\hat A_p^{(1)}(k)}{\hat B_p(k)}}\Bigg)^{-\frac{n
 \left(1-\hat D(k)\right)}{\hat A_p^{(1)}(k)}\,\frac{\hat B_p(k)}{1-\hat
 D(k)}}
\end{align}
for small $k$, for which $\hat A_p^{\sss(1)}(k)$ is bounded away from 0 and
$\hat B_p(k)$ is close to 0.  For $k_n$ defined in \refeq{kn-def},
\begin{align}\lbeq{elim}
\bigg(1+\frac{\hat B_p(k_n)}{\hat A_p^{\sss(1)}(k_n)}\bigg)^{\frac{\hat
 A_p^{\sss(1)}(k_n)}{\hat B_p(k_n)}}\!\!\underset{n\uparrow\infty}\to e,&&
\frac{n(1-\hat D(k_n))}{\hat A_p^{\sss(1)}(k_n)}\underset{n\uparrow
 \infty}\to\frac{|k|^{\alpha\wedge2}}{\hat A_p^{\sss(1)}(0)},
\end{align}
where we have used the continuity: $\hat A_p^{\sss(1)}(k_n)\to\hat A_p^{\sss
(1)}(0)=1+O(\lambda)$.  By \refeq{piNcos-bd} and \refeq{B-def}, $\hat B_p(k)
/(1-\hat D(k))=1+O(\lambda)$ uniformly in $k$.  This completes the proof of
Theorem~\ref{thm:limiting}(ii) using \refeq{limit-main}. \QED

\section{Review of the lace expansion}\label{s:laceexp}
\subsection{Derivation of the expansion}
In this section, we briefly explain the lace expansion \refeq{laceexp-x} for
oriented percolation.  In the literature, there are currently three different
ways to obtain \refeq{laceexp-x} and different representations for
$\pi_p(x,n)$.  One is based on an algebraic approach using the Markov property
\cite{NY1}, another one is to use inclusion-exclusion and nested expectations
\cite{RHS1}, and the other is to use inclusion-exclusion and the Markov
property \cite{S1}.  Here, we provide a quick overview of the third approach,
which is thought to be conceptually simplest.  The readers who are familiar to
the lace expansion for oriented percolation may skip this section and
immediately go to Section~\ref{s:bootstrapping}.

Recall that $\varphi_p(x,n)$ is the probability that $(o,0)$ is connected to
$(x,n)$.  In order for this event to occur, there are two disjoint events
depending on whether there is or is not a pivotal bond for $\{(o,0)\to(x,n)\}$.
If a bond $b$ is pivotal for $\{(o,0)\to(x,n)\}$, then $(x,n)$ is not contained
in the set of sites connected from $(o,0)$ without using $b$.  For
$(v,l)\in\Zd\times\Zp$, let
\begin{align}
\tilde\cC^b(v,l)=\{(y,n)\in\Zd\times\Zp:(v,l)\to(y,n)
 \text{ without using }b\}.
\end{align}
If there is no pivotal bond for $\{(o,0)\to(x,n)\}$, then $(o,0)=(x,n)$
or there are at least two bond-disjoint nonzero occupied paths from
$(o,0)$ to $(x,n)$.  We denote this event by $\{(o,0)\db(x,n)\}$ and define
\begin{align}\lbeq{pi0-def}
\pi_p^{\sss(0)}(x,n)=\mP_p((o,0)\db(x,n)).
\end{align}
Then, by taking the first pivotal bond $b$ (if it exists) for
$\{(o,0)\to(x,n)\}$, we obtain
\begin{align}\lbeq{laceexp-der1}
\varphi_p(x,n)=\pi_p^{\sss(0)}(x,n)+\sum_b\mP_p\big((o,0)\db b\to(x,n)
 \notin\tilde\cC^b(o,0)\big),
\end{align}
where, by denoting $b=(\bb,\tb)$, we have used the abbreviation
\begin{align}\lbeq{abbrev1}
\{(o,0)\db b\to(x,n)\}&=\{(o,0)\db\bb\}\cap\{b\to(x,n)\}\nn\\
&=\{(o,0)\db\bb\}\cap\{b\text{ is occupied}\}
 \cap\{\tb\to(x,n)\}.
\end{align}
By inclusion-exclusion in terms of the condition
$(x,n)\notin\tilde\cC^b(o,0)$, the second term in \refeq{laceexp-der1} is
\begin{align}\lbeq{incl-excl1}
\sum_b\mP_p((o,0)\db b\to(x,n))&-\sum_b\mP_p\big((o,0)\db b\to(x,n)
 \in\tilde\cC^b(o,0)\big)\nn\\
=(\pi_p^{\sss(0)}*q_p*\varphi_p)(x,n)&-R_p^{\sss(1)}(x,n)
\end{align}
where we have applied the Markov property for the first term, and
\begin{align}\lbeq{R1-def}
R_p^{\sss(1)}(x,n)=\sum_b\mP_p\big((o,0)\db b\to(x,n)\in\tilde\cC^b
 (o,0)\big).
\end{align}
Therefore,
\begin{align}
\varphi_p(x,n)=\pi_p^{\sss(0)}(x,n)+(\pi_p^{\sss(0)}*q_p*\varphi_p)(x,n)
 -R_p^{\sss(1)}(x,n).
\end{align}
This completes the first step of the full expansion \refeq{laceexp-x}.

To proceed the expansion further, it suffices to consider $R_p^{\sss(1)}(x,n)$.
Given a set $\cC$ of vertices, we define
\begin{align}
E(b,(x,n);\cC)&=\{b\to(x,n)\in\cC\}\cap\big\{\nexists b'\text{ pivotal for }
 \{\tb\to(x,n)\}\text{ satisfying }\bb'\in\cC\big\}.
\end{align}
and, for $N\ge1$ and $\vec b_N=(b_1,\dots,b_N)$,
\begin{align}
\tilde E_{\vec b_N}^{\sss(N)}(x,n)=\{(o,0)\db\bb_1\}\cap\bigcap_{i=1}^N
 E\big(b_i,\bb_{i+1};\tilde\cC^{b_i}(\tb_{i-1})\big),
\end{align}
with the convention $\tb_0=(o,0)$ and $\bb_{N+1}=(x,n)$.  For $N\ge0$, we
define
\begin{align}
\pi_p^{\sss(N)}(x,n)&=\begin{cases}
 \dpst\mP_p((o,0)\db(x,n))&(N=0),\\[7pt]
 \dpst\sum_{\vec b_N}\mP_p\big(\tilde E_{\vec b_N}^{\sss(N)}(x,n)\big)\quad
  &(N\ge1),
 \end{cases}\lbeq{piN-def}\\[3pt]
R_p^{\sss(N+1)}(x,n)&=\sum_{\vec b_{N+1}}\mP_p\Big(\tilde E_{\vec b_N}^{\sss
 (N)}(\bb_{N+1})\cap\big\{b_{N+1}\to(x,n)\in\tilde\cC^{b_{N+1}}(\tb_N)\big\}
 \Big),\lbeq{RN-def}
\end{align}
which are consistent with \refeq{pi0-def} and \refeq{R1-def}.  It has been
proved \cite{vdHS04,S1} that
\begin{align}\lbeq{laceexp-der3}
R_p^{\sss(N)}(x,n)=\pi_p^{\sss(N)}(x,n)+(\pi_p^{\sss(N)}*q_p*\varphi_p)(x,n)
 -R_p^{\sss(N+1)}(x,n).
\end{align}
We note that $R_p^{\sss(N)}(x,n)$ involves the sum over $b_1,\dots,b_N$ with
$\tb_{j-1}<\tb_j$ for $j=2,\dots,N$, hence $R_p^{\sss(N)}(x,n)=0$ if $N>n$.
Repeatedly using \refeq{laceexp-der3}, we arrive at the following conclusion:

\begin{proposition}[\cite{vdHS04,S1}]\label{prop:laceexp-2pt}
\begin{align}\lbeq{lace-prop}
\varphi_p(x,n)=\pi_p(x,n)+(\pi_p*q_p*\varphi_p)(x,n),
\end{align}
where
\begin{align}\lbeq{PiN-def}
\pi_p(x,n)=\sum_{N=0}^\infty(-1)^N\pi_p^{\sss(N)}(x,n).
\end{align}
\end{proposition}

\begin{figure}[t]
\begin{gather*}
\pi_p^{\sss(0)}(x,n):~
 \raisebox{-3.3pc}{\includegraphics[scale=0.2]{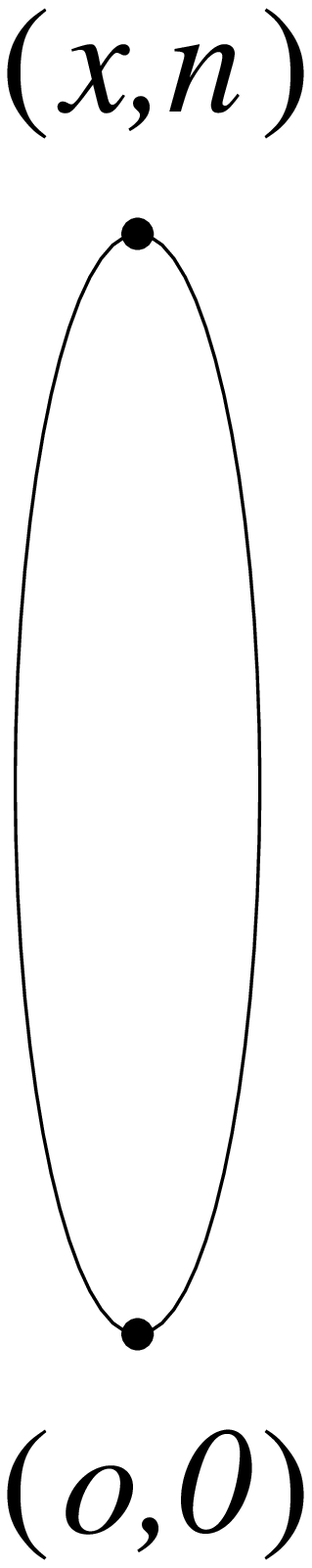}}\qquad\qquad
\pi_p^{\sss(1)}(x,n):~
 \raisebox{-3.3pc}{\includegraphics[scale=0.2]{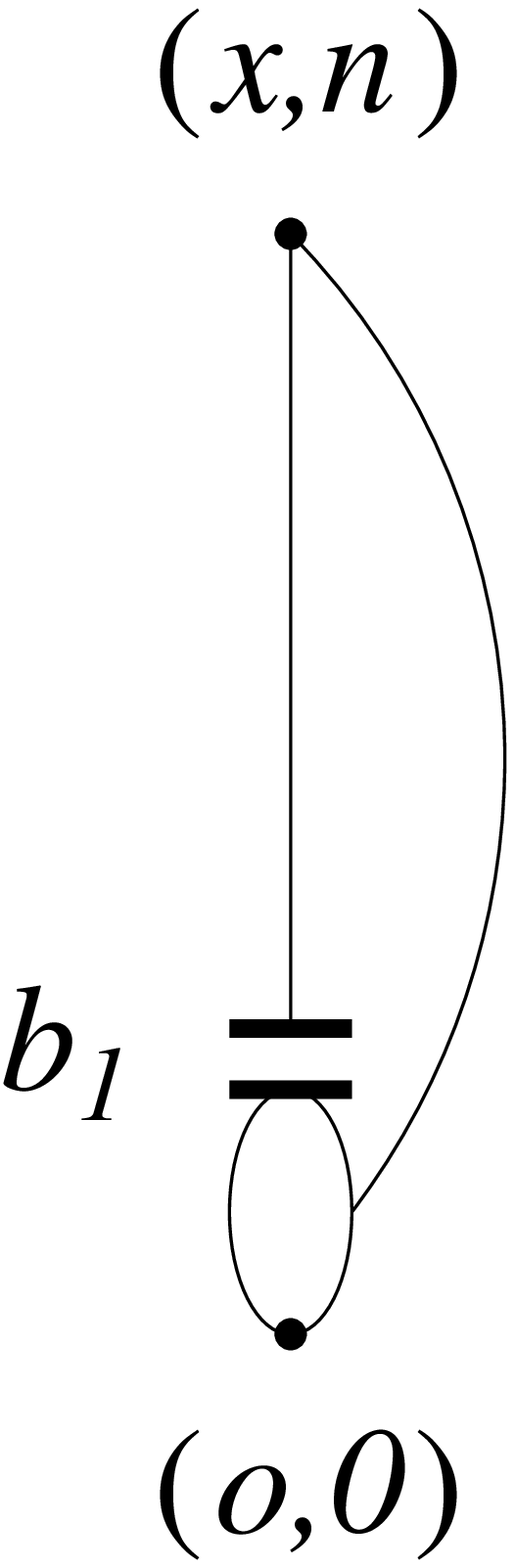}}\qquad\qquad
\pi_p^{\sss(2)}(x,n):~
 \raisebox{-3.3pc}{\includegraphics[scale=0.2]{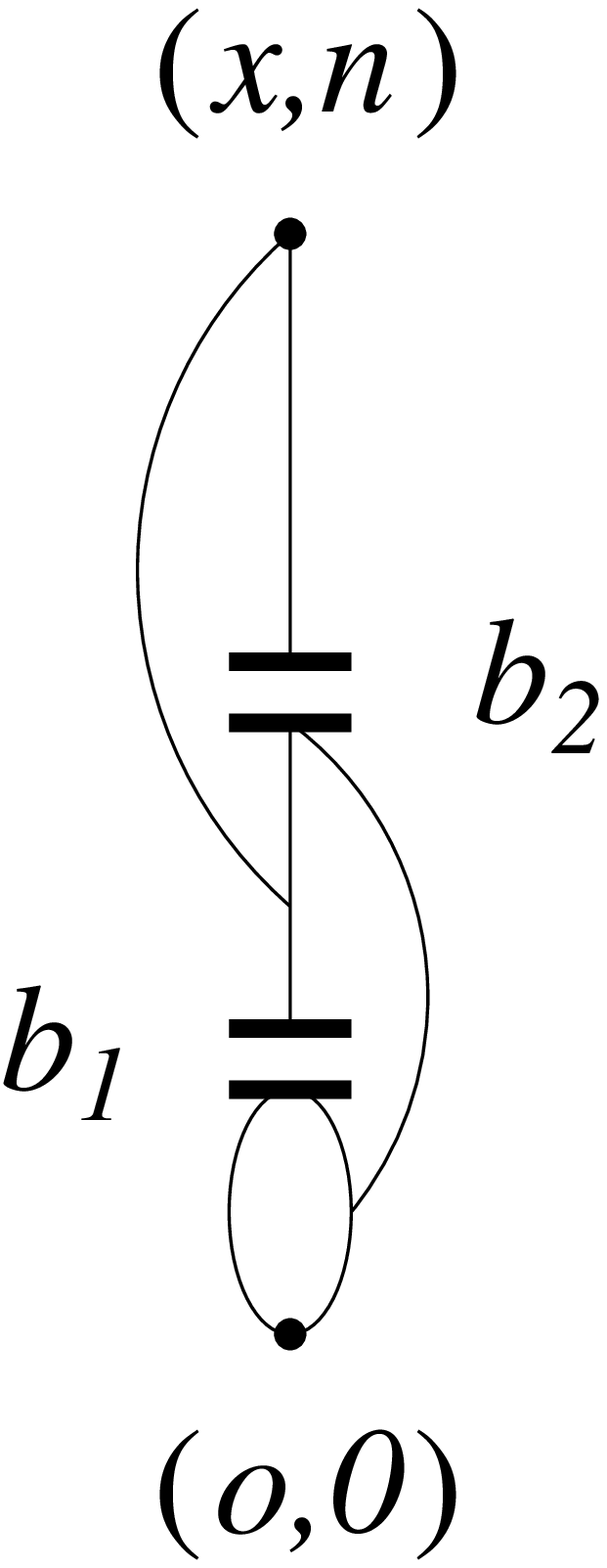}}\quad\bigcup
 \quad\raisebox{-3.3pc}{\includegraphics[scale=0.2]{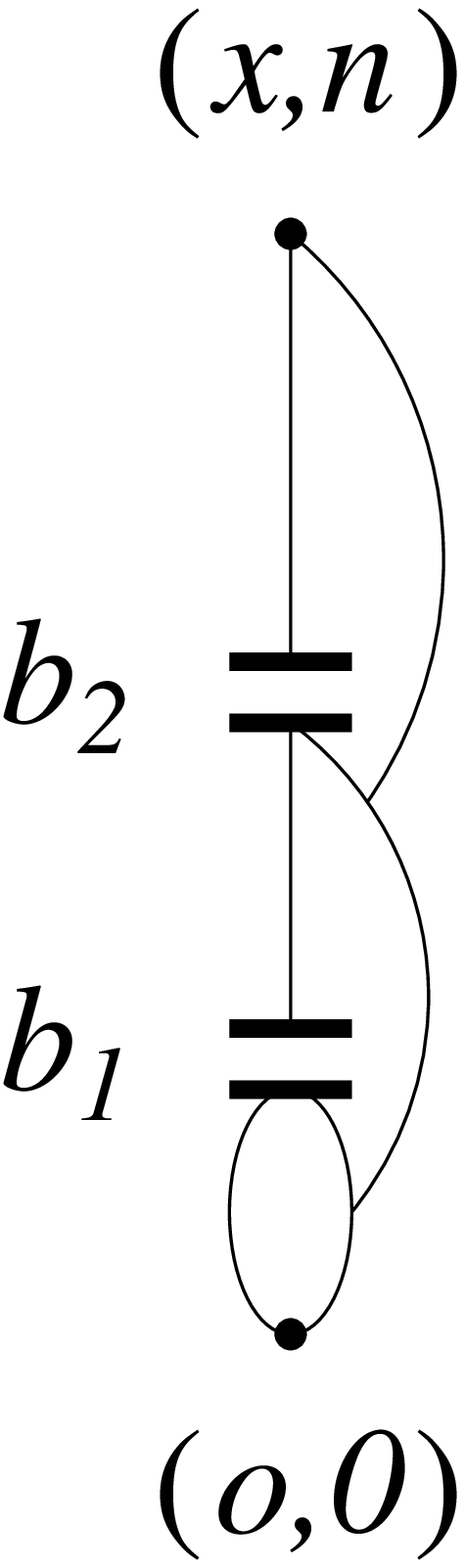}}\\
\Pi_p^{\sss(1)}(x,n):~
 \raisebox{-3.3pc}{\includegraphics[scale=0.2]{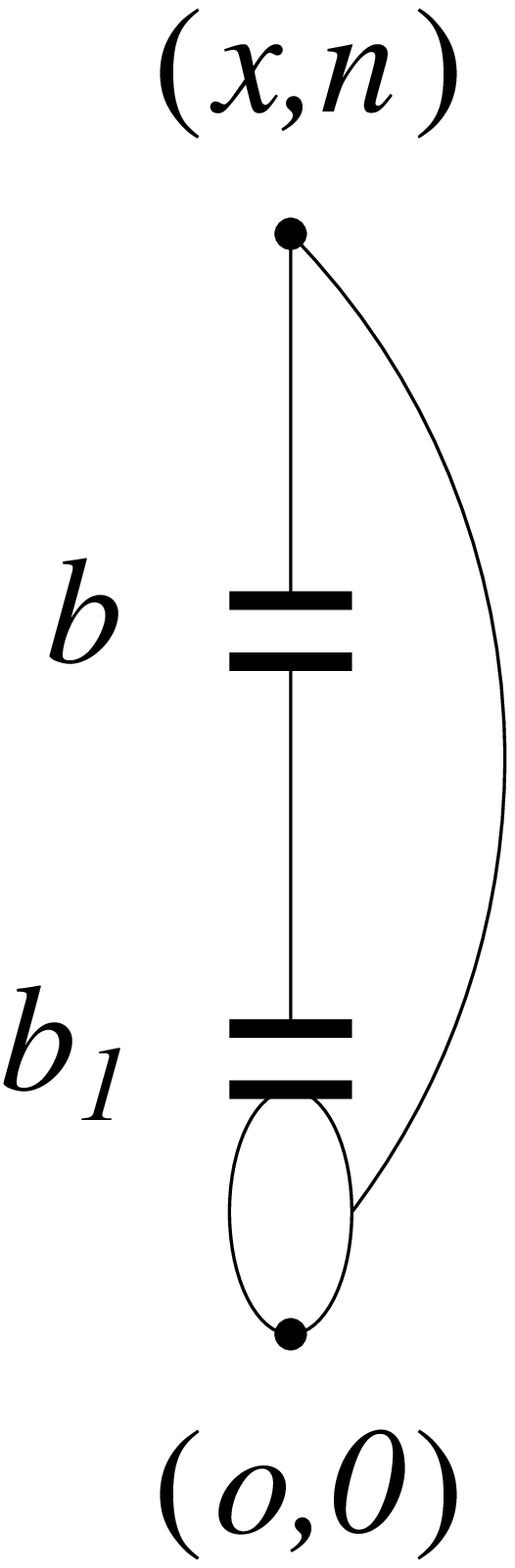}}\qquad\qquad
\Pi_p^{\sss(2)}(x,n):~
 \raisebox{-3.3pc}{\includegraphics[scale=0.2]{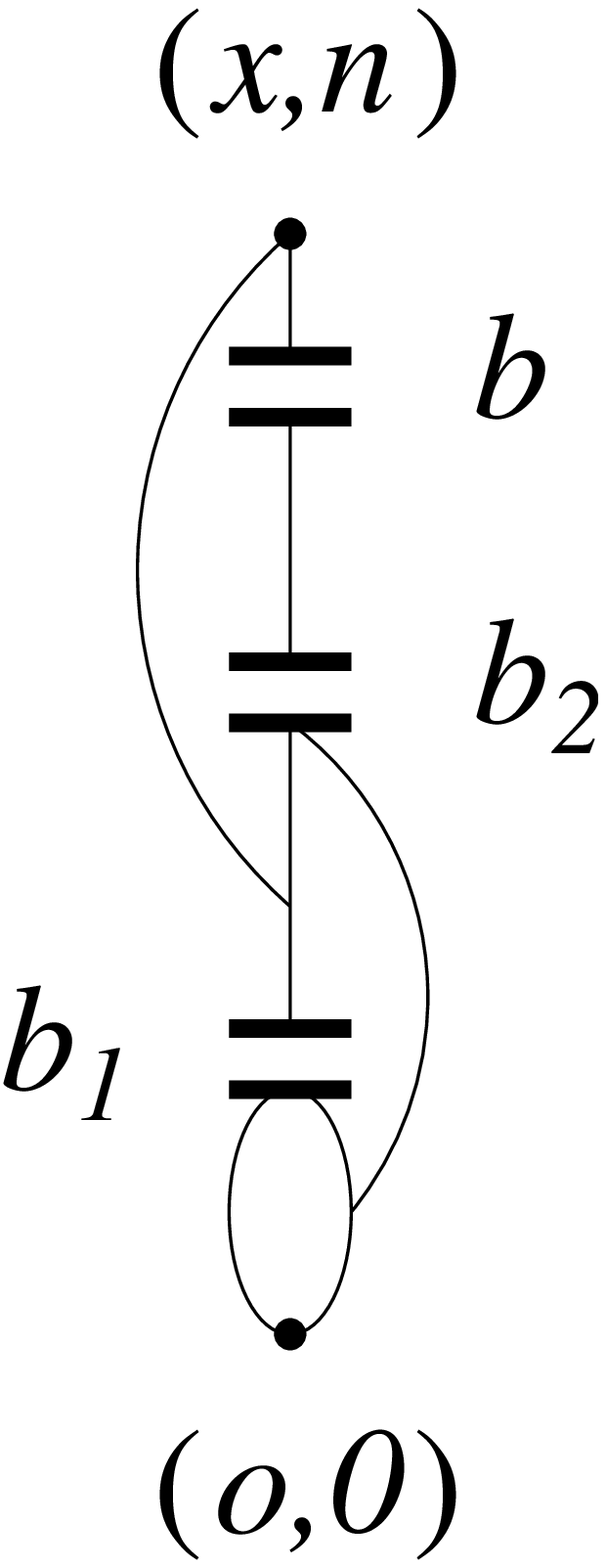}}\quad\bigcup
 \quad\raisebox{-3.3pc}{\includegraphics[scale=0.2]{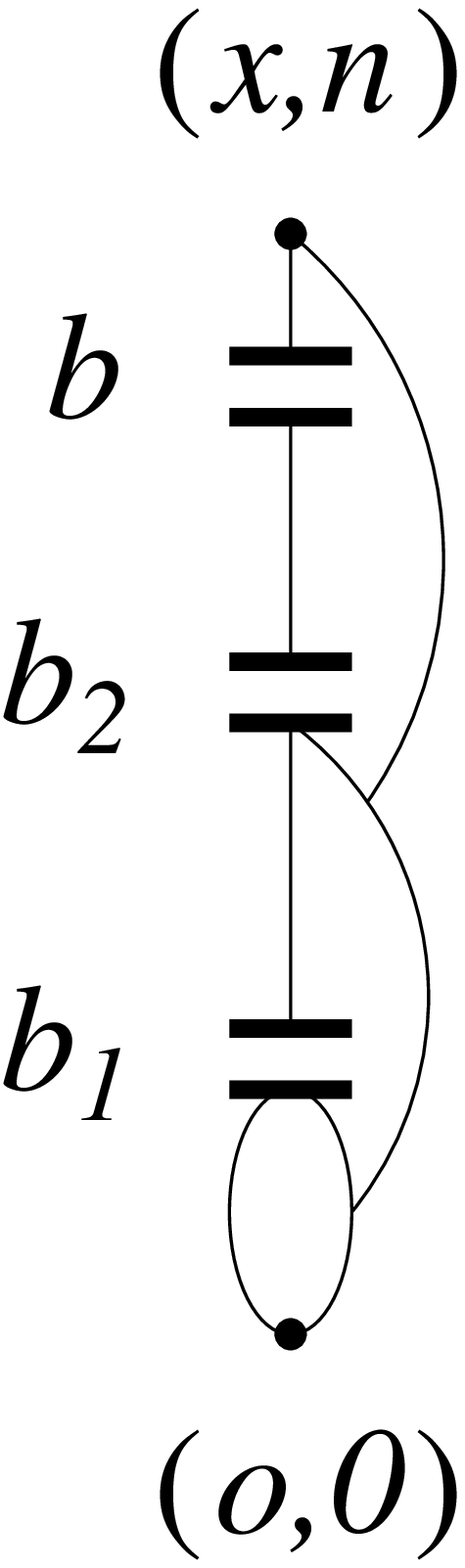}}
\end{gather*}
\caption{\label{fig-pi012}Schematic representations of $\pi_p^{\sss(N)}(x,n)$
for $N=0,1,2$ and $\Pi_p^{\sss(N)}(x,n)$ for $N=1,2$.  The $b$'s are bonds that
are summed over.}
\end{figure}

Extending the above idea, we obtain the following
representation\footnote{\label{FN:Russo}Proposition~\ref{prop:laceexp-der} is a
result of applying Russo's formula \cite{r81} to $\varphi_p(x,n)$ and compare
the result with the derivative of \refeq{lace-prop}.  Since Russo's formula can
be used only for finite systems, we should first approximate $\varphi_p(x,n)$
by a finite-volume version $\varphi_{p,R}(x,n)\equiv\mP_p((o,0)\to(x,n)$ in
$\Lambda_R)$, where $\Lambda_R=(\mZ\cap[-R,R])^d\times\Zp$, and then apply
Russo's formula. This strategy is explained in \cite[Section~3.2]{vdHS04},
where a sort of finite-confinement argument of random-walk paths is used. Since
the tail of the underlying random walk in the current setting does not decay
fast, we restrict $p$ to $p<\pc$ and use the fact that $\chi_p<\infty$ and
$\tilde\chi_{p,R}\equiv\sum_{(x,n)\notin\Lambda_R}\varphi_p(x,n)\to0$ as
$R\to\infty$.  Then, the corresponding quantities to the first and second lines
of \cite[(3.58)]{vdHS04} are bounded respectively by $\tilde\chi_{p,R}$ and
$\chi_p^3\tilde\chi_{p,R}$, both of which tend to zero as $R\to\infty$, hence
we obtain \refeq{pider-exp}--\refeq{Pider-exp}.} of $\partial_p\pi_p(x,n)$ for
$p\in(0,\pc)$, which will be used in Section~\ref{s:pc} to prove
Proposition~\ref{prop:piderbd}.

\begin{proposition}[\cite{vdHS04}]\label{prop:laceexp-der}
For $p\in(0,\pc)$,
\begin{align}\lbeq{pider-exp}
\partial_p\pi_p(x,n)=\frac1p\sum_{N=1}^\infty(-1)^N\Pi_p^{\sss(N)}(x,n),
\end{align}
where
\begin{align}\lbeq{Pider-exp}
\Pi_p^{\sss(N)}(x,n)=\sum_{\vec b_N,b}\sum_{j=1}^N\mP_p\Big(\tilde E_{\vec
 b_N}^{\sss(N)}(x,n)\cap\big\{b=b_j\text{ or $b$ is pivotal for }\{\tb_j\to
 \bb_{j+1}\}\big\}\Big),
\end{align}
with the convention $\bb_{N+1}=(x,n)$.
\end{proposition}

\subsection{Diagrammatic bounds on the expansion coefficients}\label{ss:diagbd}
In this section, we provide diagrammatic bounds on $\pi_p^{\sss(N)}(x,n)$ and
$\Pi_p^{\sss(N)}(x,n)$.  These bounds consist of two-point functions, and are
results of applications of the BK inequality~\cite{BK} and
\begin{align}
\varphi_p(x,n)\le(q_p*\varphi_p)(x,n)\qquad(n\ge1).
\end{align}
For example, $\pi_p^{\sss(0)}(x,n)$ is bounded as
\begin{align}\lbeq{pi0-BKappl}
\pi_p^{\sss(0)}(x,n)\le\varphi_p(x,n)^2=\delta_{x,o}\delta_{n,0}+\big((1-
 \delta_{x,o}\delta_{n,0})\,\varphi_p(x,n)\big)^2\le\delta_{x,o}\delta_{n,0}
 +(q_p*\varphi_p)(x,n)^2.
\end{align}
The other terms are bounded similarly.

Let $\varphi_p^{\sss(m)}(x,n)=\varphi_p(x,n)m^n$ and define the weighted bubble
$W_p^{\sss(m)}(k)$, the triangles $T_p^{\sss(m)}$ and $\tilde T_p$, the square
$S_p^{\sss(m)}$ and the H-shaped diagrams $H_p$ as (see
Figure~\ref{fig-tildeT&H})
\begin{align}
W_p^{\sss(m)}(k)&=\sup_{(x,n)}\sum_{(y,t)}\big(1-\cos(k\cdot y)\big)\times
 \begin{cases}
 (q_p*\varphi_p)(y,t)\cdot(mq_p*\varphi_p^{\sss(m)})(y-x,t-n),
  &\text{if }m<1,\\
 (mq_p*\varphi_p^{\sss(m)})(y,t)\cdot(q_p*\varphi_p)(y-x,t-n),
  &\text{if }m\ge1,
 \end{cases}\\
T_p^{\sss(m)}&=\sup_{(x,n)}\sum_{(y,t)}(q_p*\varphi_p*\varphi_p)(y,t)
 \cdot(mq_p*\varphi_p^{\sss(m)})(y-x,t-n),\\
S_p^{\sss(m)}&=\sup_{(x,n)}\sum_{(y,t)}(q_p*\varphi_p*\varphi_p*
 \varphi_p)(y,t)\cdot(mq_p*\varphi_p^{\sss(m)})(y-x,t-n),\lbeq{S-def}\\
\tilde T_p&=\sup_{(x,n)}\sum_{(y,t)}(q_p*\varphi_p*q_p*\varphi_p)(y,t)\cdot
 (q_p*\varphi_p)(y-x,t-n),\lbeq{tildeT-def}\\
H_p&=\sup_{(x,n),(x',n')}\sum_{(y_i,t_i),\,i=1,2,3}(q_p*\varphi_p)(y_1,t_1)
 \cdot(\varphi_p*q_p*\varphi_p)(y_2-y_1,t_2-t_1)\nn\\
&\hspace{9pc}\times(q_p*\varphi_p)(y_2-x,t_2-n)\cdot(q_p*\varphi_p)(y_3-y_1,
 t_3-t_1)\nn\\
&\hspace{9pc}\times(q_p*\varphi_p)(x'+y_3-y_2,n'+t_3-t_2).\lbeq{H-def}
\end{align}
The expansion coefficients obey the following bounds:
\begin{figure}[t]
\begin{align*}
\tilde T_p=\sup_{(x,n)\in\mZ^{d+1}}\quad\raisebox{-4pc}{\includegraphics
 [scale=0.22]{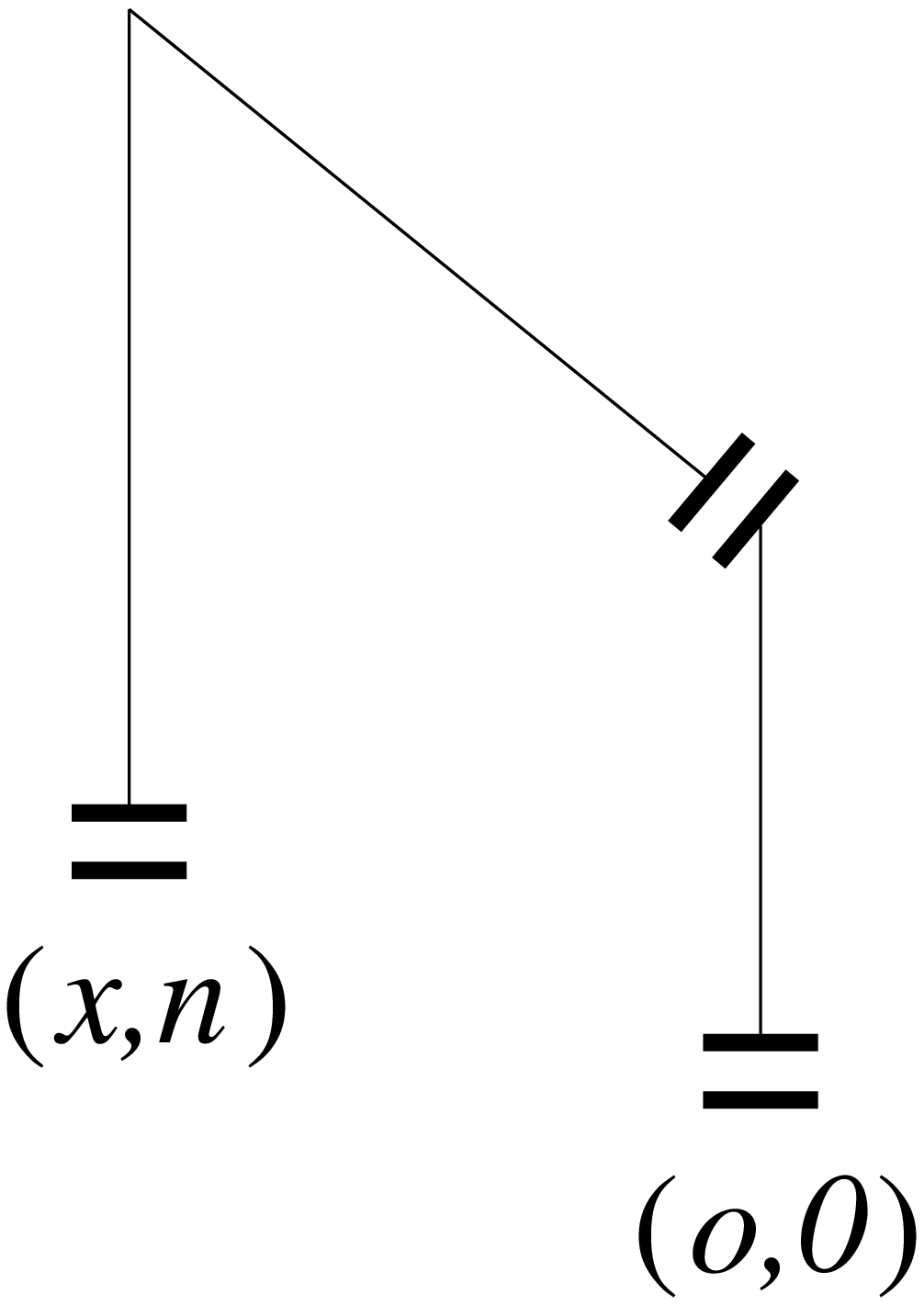}}&&&&
H_p=\sup_{(x,n),(x',n')\in\mZ^{d+1}}\quad\raisebox{-5pc}{\includegraphics
 [scale=0.22]{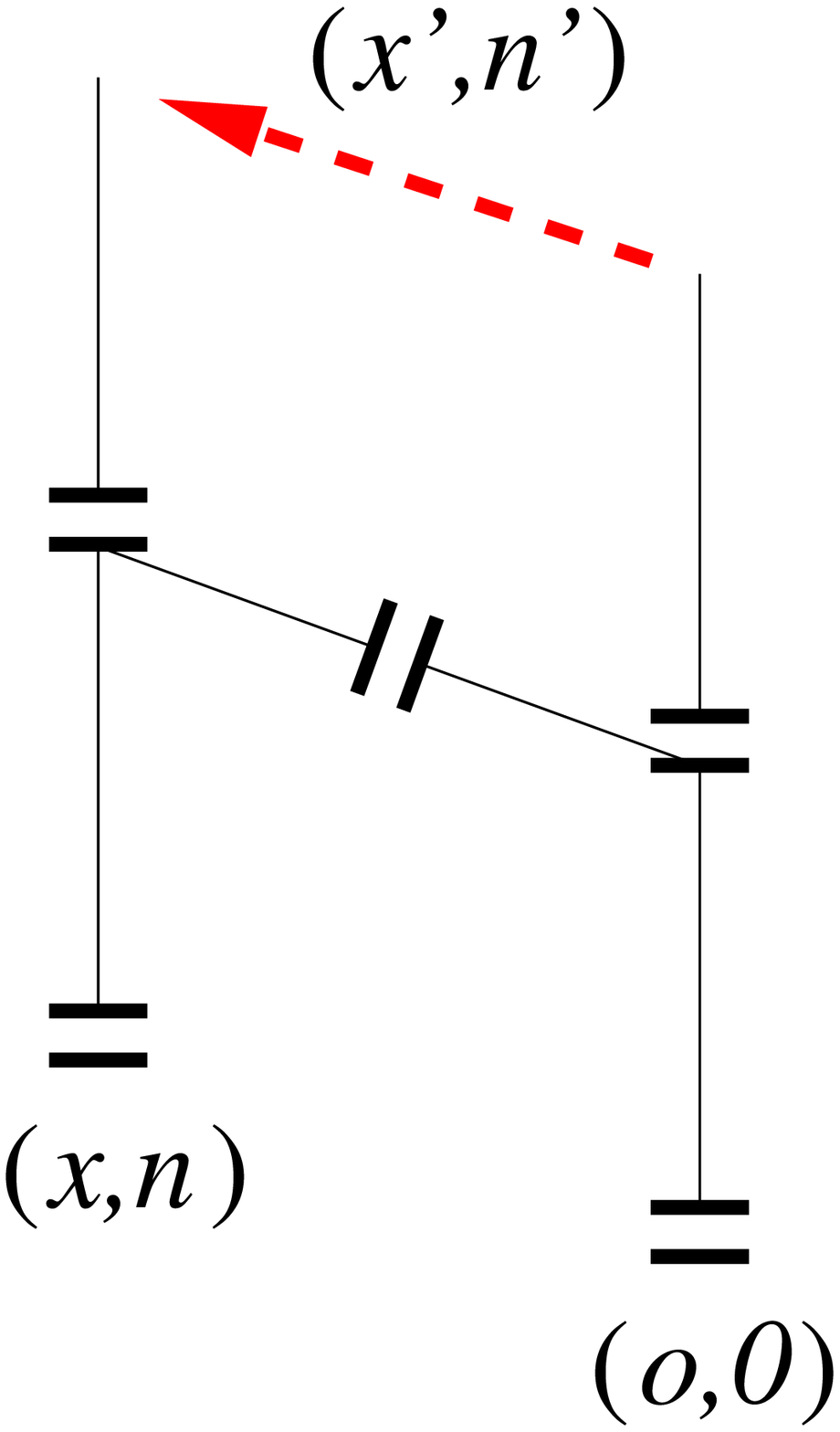}}
\end{align*}
\caption{\label{fig-tildeT&H}Schematic representations of $\tilde T_p$ and
$H_p$.}
\end{figure}

\begin{proposition}\label{prop:diagbd}
\begin{enumerate}[(i)]
\item
For $N\ge0$ and $r=0,1,2$,
\begin{align}
&\sum_{(x,n)\in\Zd\times\mN}n^r\pi_p^{\sss(N)}(x,n)m^n\le(N+1)^r(1+2T_p^{\sss
 (m)})(2T_p^{\sss(m)})^{(N-1)\vee0}\times
 \begin{cases}
 T_p^{\sss(m)}&(r=0,1),\\
 S_p^{\sss(m)}&(r=2),
 \end{cases}\lbeq{piNn-diagbd}\\
&\sum_{(x,n)\in\Zd\times\Zp}\big(1-\cos(k\cdot x)\big)\pi_p^{\sss(N)}(x,n)m^n
 \le3(N+1)^2(1+2T_p^{\sss(m)})(2T_p^{\sss(m)})^{(N-1)\vee0}W_p^{\sss(m)}(k).
 \lbeq{piNcos-diagbd}
\end{align}
\item For $N\ge1$,
\begin{align}\lbeq{PiN-diagbd}
\sum_{(x,n)\in\Zd\times\Zp}\Pi_p^{\sss(N)}(x,n)\le N(1+2T_p^{\sss(1)})\Big(
 (T_p^{\sss(1)}+\tilde T_p)(2T_p^{\sss(1)})^{N-1}+H_p(2T_p^{\sss(1)})^{(N-2)
 \vee0}\Big).
\end{align}
\end{enumerate}
\end{proposition}

The proof of the above proposition is irrelevant in this paper, and is found in
\cite{s??}.

\section{Proof of Proposition~\ref{prop:bootstrapping}}\label{s:bootstrapping}
In this section, we prove Proposition~\ref{prop:bootstrapping} that was the key
for the proof of Theorem~\ref{thm:IRbound}.  First, in
Section~\ref{ss:IR(iii)}, we prove Proposition~\ref{prop:bootstrapping}(iii)
that is nothing to do with the lace expansion.  Then, in
Section~\ref{ss:IR(ii)}, we prove Proposition~\ref{prop:bootstrapping}(ii)
using the trigonometric technique in \cite[Section~5.1]{Slad06}.  Finally, in
Section~\ref{ss:IR(i)}, we prove Proposition~\ref{prop:bootstrapping}(i) using
the diagrammatic bounds on the expansion coefficients in
Section~\ref{ss:diagbd}.

\subsection{Proof of Proposition~\ref{prop:bootstrapping}(iii)}\label{ss:IR(iii)}
First we prove $f(0,1)=1$.  When $p=0$, by definition we have $f_1(0,1)=0$,
$\hat\varphi_0(k,z)\equiv1$, $\mu_0(z)\equiv0$ (cf., \refeq{mu-def}) and hence
$\hat G_{\mu_0(z)}(k)\equiv1$.  Therefore, $f_2(0,1)=1$ and $f_3(0,1)=0$.

Next we discuss the continuity of $f(p,m)$.  Since $f_1(p,m)\equiv p(m\vee1)$
is obviously continuous in $p$ and $m$, we only need to investigate $f_2(p,m)$
and $f_3(p,m)$.

Fix $p<\pc$.  To prove the continuity of $f(p,m)$ in $m<m_p$, it suffices to
show that $f(p,m)$ is continuous in $m\in[0,\tilde m]$ for every $\tilde
m<m_p$.  To prove this for $f_2(p,m)$, it suffices to show that the derivative
\begin{align}\lbeq{f2-derm}
\partial_m\frac{\hat\varphi_p(k,me^{i\theta})}{\hat G_{\mu_p(me^{i\theta})}
 (k)}=\frac{\partial_m\hat\varphi_p(k,me^{i\theta})}{\hat G_{\mu_p(me^{i
 \theta})}(k)}-\hat\varphi_p(k,me^{i\theta})\frac{\partial_m\hat G_{\mu_p(m
 e^{i\theta})}(k)}{\hat G_{\mu_p(me^{i\theta})}(k)^2}
\end{align}
is bounded uniformly in $(k,\theta)\in[-\pi,\pi]^{d+1}$ and $m\in[0,\tilde m]$
(cf., \cite[Lemma~5.13]{Slad06}).  However, by
$n\varphi_p(x,n)\le(q_p*\varphi_p*\varphi_p)(x,n)$ (cf., \cite[(5.17)]{S1}), we
have
\begin{align}\lbeq{varphi-derm}
|\partial_m\hat\varphi_p(k,me^{i\theta})|\le\sum_{(x,n)}n\varphi_p(x,n)m^{n-1}
 \le p\hat\varphi_p(0,m)^2\le p\hat\varphi_p(0,\tilde m)^2.
\end{align}
Since $|\hat G_{\mu_p(me^{i\theta})}(k)|\ge\frac12$, the first term on the
right-hand side of \refeq{f2-derm} is indeed uniformly bounded.  Also, since
$\hat\varphi_p(0,m)~(\ge1)$ is nondecreasing in $m$, we obtain
\begin{align}
\bigg|\frac{\partial_m\hat G_{\mu_p(me^{i\theta})}(k)}{\hat G_{\mu_p(m
 e^{i\theta})}(k)^2}\bigg|=|\hat D(k)\,\partial_m\mu_p(me^{i\theta})|
 \le\frac{\partial_m\hat\varphi_p(0,m)}{\hat\varphi_p(0,m)^2},
\end{align}
which is uniformly bounded by $p$, as described in \refeq{varphi-derm}.
Consequently, \refeq{f2-derm} is uniformly bounded by $p\hat\varphi_p(0,\tilde
m)(2\hat\varphi_p(0,\tilde m)+1)$.  This completes the proof of the continuity
of $f_2(p,m)$ in $m\in[0,\tilde m]$.

Similarly to the above, we can easily show that the derivative
\begin{align}\lbeq{f3-derm}
\partial_m\frac{\hat G_{\mu_p(m)}(k)~\big(\hat\varphi_p(l,me^{i\theta})-\frac
 12(\hat\varphi_p(l+k,me^{i\theta})+\hat\varphi_p(l-k,me^{i\theta}))\big)}
 {\hat G_{\mu_p(me^{i\theta})}(l+jk)\,\hat G_{\mu_p(me^{i\theta})}(l+j'k)}
\end{align}
is bounded uniformly in $(k,\theta)\in[-\pi,\pi]^{d+1}$,
$(j,j')=(0,\pm1),(1,-1)$ and $m\in[0,\tilde m]$.  This justifies the continuity
of $f_3(p,m)$ in $m\in[0,\tilde m]$.

To prove the continuity of $f(p,1)$ in $p<\pc$, it suffices to show that
$f(p,1)$ is continuous in $p\in[0,\tilde p]$ for every $\tilde p<\pc$.  First
we note that, by Russo's formula \cite{r81} (see also Footnote~\ref{FN:Russo})
and the fact that $\chi_p\equiv\hat\varphi_p(0,1)~(\ge1)$ is nondecreasing in
$p$, we have, for $|z|=1$,
\begin{gather}
|\partial_p\hat\varphi_p(k,z)|\le\sum_{(x,n)}\partial_p\varphi_p(x,n)\le
 \sum_{(x,n)}(\varphi_p*q_1*\varphi_p)(x,n)\le\chi_p^2,\\
\bigg|\frac{\partial_p\hat G_{\mu_p(z)}(k)}{\hat G_{\mu_p(z)}(k)^2}\bigg|=
 |\hat D(k)\,\partial_p\mu_p(z)|\le\frac{\partial_p\chi_p}{\chi_p^2}\le1.
\end{gather}
Since $|\hat G_{\mu_p(z)}(k)|\ge\frac12$, we obtain
\begin{align}
\bigg|\partial_p\frac{\hat\varphi_p(k,z)}{\hat G_{\mu_p(z)}(k)}\bigg|\le
 \bigg|\frac{\partial_p\hat\varphi_p(k,z)}{\hat G_{\mu_p(z)}(k)}\bigg|+
 |\hat\varphi_p(k,z)|\,\bigg|\frac{\partial_p\hat G_{\mu_p(z)}(k)}{\hat
 G_{\mu_p(z)}(k)^2}\bigg|\le\chi_{\tilde p}(2\chi_{\tilde p}+1),
\end{align}
uniformly in $k\in[-\pi,\pi]^d$, $|z|=1$ and $p\in[0,\tilde p]$.  This implies
the continuity of $f_2(p,1)$ in $p\in[0,\tilde p]$ for every $\tilde p<\pc$.

The continuity of $f_3(p,1)$ can be proved in a similar way.  This completes
the proof of Proposition~\ref{prop:bootstrapping}(iii). \QED

\subsection{Proof of Proposition~\ref{prop:bootstrapping}(ii)}\label{ss:IR(ii)}
In this section, we prove that, for every $p<\pc$ and $m<m_p$, the weaker bound
$f(p,m)\le3$ and \refeq{piN-bd}--\refeq{piNcos-bd} imply the stronger bound
$f(p,m)\le2$ when $d>2(\alpha\wedge2)$ and $L\gg1$.

First, by \refeq{pmppi} (recall that this is a consequence of the assumed bound
\refeq{piN-bd} and the fact that $\hat\varphi_p(0,m)$ diverges as $m\uparrow
m_p$) and \refeq{piN-bd}, we immediately obtain
\begin{align}
f_1(p,m)\equiv p(m\vee1)\le pm_p=\hat\pi_p(0,m_p)^{-1}=1+O(\lambda)\le2.
\end{align}

Next we consider $f_2(p,m)$.  First we rewrite $\hat\varphi_p(k,z)/\hat
G_{\mu_p(z)}(k)$ as
\begin{align}
\frac{\hat\varphi_p(k,z)}{\hat G_{\mu_p(z)}(k)}&=\hat\pi_p(k,z)+\hat
 \varphi_p(k,z)\bigg(\frac1{\hat G_{\mu_p(z)}(k)}-\frac{\hat\pi_p(k,z)}
 {\hat\varphi_p(k,z)}\bigg)\nn\\
&=\hat\pi_p(k,z)+\hat\varphi_p(k,z)\big(pz\hat\pi_p(k,z)-\mu_p(z)\big)
 \hat D(k)\nn\\
&=\hat\pi_p(k,z)+\hat\varphi_p(k,z)\bigg(p|z|\hat\pi_p(k,z)-1+\frac1
 {\hat\varphi_p(0,|z|)}\bigg)e^{i\arg(z)}\hat D(k),
\end{align}
where
\begin{align}
p|z|\hat\pi_p(k,z)-1+\frac1{\hat\varphi_p(0,|z|)}&=p|z|\big(\hat\pi_p(k,
 z)-\hat\pi_p(0,|z|)\big)-\big(\underbrace{1-p|z|\hat\pi_p(0,|z|)}_{\hat
 \pi_p(0,|z|)/\hat\varphi_p(0,|z|)}\big)+\frac1{\hat\varphi_p(0,|z|)}\nn\\
&=p|z|\big(\pi_p(k,z)-\pi_p(0,|z|)\big)+\frac{1-\hat\pi_p(0,|z|)}
 {\hat\varphi_p(0,|z|)}.
\end{align}
We note that $|\hat\pi_p(k,z)-1|=O(\lambda)$, due to \refeq{piN-bd} for $r=0$,
and that $|\hat\varphi_p(k,z)/\hat\varphi_p(0,|z|)|\le1$ by definition.  To
complete the proof of $f_2(p,m)=1+O(\lambda)\le2$, it thus suffices to show
that
\begin{align}\lbeq{f2-(ii)}
|\hat\varphi_p(k,z)|\Big(|\pi_p(k,z)-\pi_p(0,z)|+|\pi_p(0,z)-\pi_p(0,|z|)|
 \Big)=O(\lambda),
\end{align}
uniformly in $k\in[-\pi,\pi]^d$ and $z\in\mC$ with $|z|=m$ or 1.  However, by
\refeq{piN-bd}--\refeq{piNcos-bd} and denoting $\theta=\arg(z)$, we have
\begin{align}
|\pi_p(k,z)-\pi_p(0,z)|&\le O(\lambda)\,\hat G_{\mu_p(m\vee1)}(k)^{-1}\le
 O(\lambda)\big(1-\mu_p(m\vee1)+1-\hat D(k)\big),\\
|\pi_p(0,z)-\pi_p(0,|z|)|&=\bigg|\sum_{(x,n)}\pi_p(x,n)|z|^n(e^{i\theta n}
 -1)\bigg|\le|\theta|\sum_{(x,n)}n|\pi_p(x,n)||z|^n=O(\lambda)|\theta|.
\end{align}
On the other hand, by $f_2(p,m)\le3$, \refeq{G-IRbd} and
$|\mu_p(z)|\le\mu_p(m\vee1)$ for $|z|=m$ or 1 (cf., \refeq{mu-def}),
\begin{align}
|\hat\varphi_p(k,z)|\le\frac{3c}{1-\mu_p(m\vee1)+|\theta|+1-\hat D(k)}.
\end{align}
This completes the proof of \refeq{f2-(ii)}, and hence $f_2(p,m)\le2$.

For $f_3(p,m)$, we introduce the following notation for $\hat
f(l)\equiv\sum_{x\in\Zd}f(x)e^{il\cdot x}$:
\begin{align}\lbeq{laplacian}
\Delta_k\hat f(l)=\hat f(l+k)+\hat f(l-k)-2\hat f(l).
\end{align}
We note that $-\frac12\Delta_k\hat f(l)$ is the Fourier transform of
$(1-\cos(k\cdot x))f(x)$:
\begin{align}\lbeq{DeltaFourier}
-\frac12\Delta_k\hat f(l)=\sum_{x\in\Zd}f(x)\bigg(e^{il\cdot x}-\frac{
 e^{i(l+k)\cdot x}+e^{i(l-k)\cdot x}}2\bigg)=\sum_{x\in\Zd}f(x)\big(1-
 \cos(k\cdot x)\big)e^{il\cdot x}.
\end{align}
Recall the definition of $f_3(p,m)$ whose numerator contains
$-\frac12\Delta_k\hat\varphi_p(l,z)$.  Let
\begin{align}
\hat a_p(l,z)=pz\hat D(l)\,\hat\pi_p(l,z)\equiv\sum_{(x,n)}(q_p*\pi_p)
 (x,n)z^n\cos(l\cdot x),
\end{align}
so that $\hat\varphi_p(l,z)=\hat\pi_p(l,z)/(1-\hat a_p(l,z))$.  Then, we have
\begin{align}\lbeq{f3-num}
\Delta_k\hat\varphi_p(l,z)&=\frac{\Delta_k\hat\pi_p(l,z)}{1-\hat a_p(l,z)}
 +\sum_{j=\pm1}\frac{(\hat\pi_p(l+jk,z)-\hat\pi_p(l,z))(\hat a_p(l+jk,z)
 -\hat a_p(l,z))}{(1-\hat a_p(l,z))(1-\hat a_p(l+jk,z))}\nn\\
&\quad+\hat\pi_p(l,z)\,\Delta_k\frac1{1-\hat a_p(l,z)},
\end{align}
where, by \refeq{piN-bd}--\refeq{piNcos-bd} and $f_2(p,m)\le2$,
\begin{align}\lbeq{f3-num-2nd-1}
\bigg|\frac{\Delta_k\hat\pi_p(l,z)}{1-\hat a_p(l,z)}\bigg|=\bigg|\frac{
 \Delta_k \hat\pi_p(l,z)}{\hat\pi_p(l,z)}\bigg||\hat\varphi_p(l,z)|\le
 O(\lambda)\,\hat G_{\mu_p(m\vee1)}(k)^{-1}|\hat G_{\mu_p(z)}(l)|.
\end{align}

The second term of \refeq{f3-num} can be bounded as follows.  First, by
$|e^{il\cdot x}(e^{ijk\cdot x}-1)|\le|\sin(k\cdot x)|+1-\cos(k\cdot x)$ for
$j=\pm1$,
\begin{align}
|\hat\pi_p(l+jk,z)-\hat\pi_p(l,z)|\le\sum_{(x,n)}|\sin(k\cdot x)||\pi_p(x,n)|
 |z|^n+\sum_{(x,n)}\big(1-\cos(k\cdot x)\big)|\pi_p(x,n)||z|^n,
\end{align}
where the second term is bounded by $O(\lambda)\hat G_{\mu_p(m\vee1)}(k)^{-1}$,
due to \refeq{piNcos-bd}.  By the Cauchy-Schwarz inequality and using
\refeq{piN-bd}--\refeq{piNcos-bd}, the first term is bounded by
\begin{align}
&\bigg(\sum_{(x,n):x\ne o}|\pi_p(x,n)||z|^n\bigg)^{1/2}\bigg(\sum_{(x,n):x\ne
 o}\sin^2(k\cdot x)|\pi_p(x,n)||z|^n\bigg)^{1/2}\nn\\
&\quad\le O(\lambda)^{1/2}\bigg(\sum_{(x,n)}\big(1-\cos(k\cdot x)\big)|\pi_p
 (x,n)||z|^n\bigg)^{1/2}\le O(\lambda)\,\hat G_{\mu_p(m\vee1)}(k)^{-1/2}.
\end{align}
Therefore, $|\hat\pi_p(l+jk,z)-\hat\pi_p(l,z)|\le O(\lambda)\hat
G_{\mu_p(m\vee1)}(k)^{-1/2}$.  Similarly, we can show $|\hat a_p(l+jk,z)-\hat
a_p(l,z)|\le O(1)\hat G_{\mu_p(m\vee1)}(k)^{-1/2}$, where we use
\begin{align}\lbeq{key}
&\sum_{(x,n)}\big(1-\cos(k\cdot x)\big)(q_p*|\pi_p|)(x,n)|z|^n\nn\\
&\le5p|z|\bigg(\underbrace{\sum_y\big(1-\cos(k\cdot y)\big)D(y)}_{1-\hat D(k)}
 \;\underbrace{\sum_{(x,n)}|\pi_p(x-y,n-1)||z|^{n-1}}_{1+O(\lambda)}\nn\\
&\quad\qquad+\sum_yD(y)\underbrace{\sum_{(x,n)}\Big(1-\cos\big(k\cdot(x-y)\big)
 \Big)|\pi_p(x-y,n-1)||z|^{n-1}}_{O(\lambda)\,\hat G_{\mu_p(m\vee1)}(k)^{-1}}
 \bigg)\nn\\
&\le10\big(2+O(\lambda)\big)\,\hat G_{\mu_p(m\vee1)}(k)^{-1}.
\end{align}
Here, the first inequality is due to $1-\cos(X+Y)\le5(1-\cos X)+5(1-\cos Y)$
(cf., \cite[(4.50)]{Slad06}), and the second inequality is due to
$f_1(p,m)\le2$ and $1-\hat D(k)\le2\hat G_{\mu_p(m\vee1)}(k)^{-1}$ (since
$\mu_p(m\vee1)\in[0,1]$).  Therefore, for $j=\pm1$,
\begin{align}\lbeq{f3-num-2nd-2}
\bigg|\frac{(\hat\pi_p(l+jk,z)-\hat\pi_p(l,z))(\hat a_p(l+jk,z)-\hat a_p(l,z))}
 {(1-\hat a_p(l,z))(1-\hat a_p(l+jk,z))}\bigg|\le O(\lambda)\,\hat G_{\mu_p(m
 \vee1)}(k)^{-1}|\hat G_{\mu_p(z)}(l)\,\hat G_{\mu_p(z)}(l+jk)|.
\end{align}

To complete bounding $\Delta_k\hat\varphi_p(l,z)$, it remains to investigate
$\Delta_k(1-\hat a_p(l,z))^{-1}$ in the last term of \refeq{f3-num}.  Let
\begin{align}
\hat a_p^{\cos}(l,z;k)&=\sum_{(x,n)}(q_p*\pi_p)(x,n)z^n\cos(l\cdot x)\cos
 (k\cdot x),\\
\hat a_p^{\sin}(l,z;k)&=\sum_{(x,n)}(q_p*\pi_p)(x,n)z^n\sin(l\cdot x)\sin
 (k\cdot x).
\end{align}
Then, by \cite[Lemma~5.3]{bchss05},
\begin{align}\lbeq{Delta1-a}
\Delta_k\frac1{1-\hat a_p(l,z)}=\frac{\hat\varphi_p(l,z)}{\hat\pi_p(l,z)}
 \bigg(\sum_{j=\pm1}\frac{\hat\varphi_p(l+jk,z)}{\hat\pi_p(l+jk,z)}\,
 \big(\hat a_p^{\cos}(l,z;k)-\hat a_p(l,z)\big)\nn\\
+2\prod_{j=\pm1}\frac{\hat\varphi_p(l+jk,z)}{\hat\pi_p(l+jk,z)}\,\hat
 a_p^{\sin}(l,z;k)^2\bigg),
\end{align}
where, by \refeq{key},
\begin{align}
|\hat a_p^{\cos}(l,z;k)-\hat a_p(l,z)|\le\sum_{(x,n)}\big(1-\cos(k\cdot
 x)\big)(q_p*|\pi_p|)(x,n)|z|^n\le10\big(2+O(\lambda)\big)\,\hat G_{\mu_p
 (m\vee1)}(k)^{-1}.
\end{align}
Moreover, by the Cauchy-Schwarz inequality,
\begin{align}
\hat a_p^{\sin}(l,z;k)^2&\le\bigg(\sum_{(x,n)}(q_p*|\pi_p|)(x,n)|z|^n\sin^2
 (l\cdot x)\bigg)\sum_{(x,n)}(q_p*|\pi_p|)(x,n)|z|^n\sin^2(k\cdot x)\nn\\
&\le2^2\bigg(\sum_{(x,n)}\big(1-\cos(l\cdot x)\big)(q_p*|\pi_p|)(x,n)|z|^n
 \bigg)\sum_{(x,n)}\big(1-\cos(k\cdot x)\big)(q_p*|\pi_p|)(x,n)|z|^n\nn\\
&\le20^2\big(2+O(\lambda)\big)^2\hat G_{\mu_p(m\vee1)}(l)^{-1}\hat G_{\mu_p
 (m\vee1)}(k)^{-1}.
\end{align}
As a result, since $f_2(p,m)\le2$ and $|\hat G_{\mu_p(z)}(l)|\le\hat
G_{\mu_p(m\vee1)}(l)$ for $|z|=m$ or 1, we obtain
\begin{align}\lbeq{f3-num-3rd}
\bigg|\Delta_k\frac1{1-\hat a_p(l,z)}\bigg|&\le\hat G_{\mu_p(m\vee1)}
 (k)^{-1}\bigg(40\big(2+O(\lambda)\big)\sum_{j=\pm1}|\hat G_{\mu_p(z)}
 (l)\,\hat G_{\mu_p(z)}(l+jk)|\nn\\
&\hspace{8pc}+80^2\big(2+O(\lambda)\big)^2|\hat G_{\mu_p(z)}(l+k)\,\hat
 G_{\mu_p(z)}(l-k)|\bigg)\nn\\
&\le2K\big(1+O(\lambda)\big)\,\hat G_{\mu_p(m\vee1)}(k)^{-1}\sum_{(j,
 j')=(0,\pm1),(1,-1)}|\hat G_{\mu_p(z)}(l+jk)\,\hat G_{\mu_p(z)}(l+j'k)|,
\end{align}
where $K=2\cdot80^2$.

Finally, by summarizing \refeq{f3-num}--\refeq{f3-num-2nd-1},
\refeq{f3-num-2nd-2} and \refeq{f3-num-3rd}, we arrive at
\begin{align}
\frac{\hat G_{\mu_p(m\vee1)}(k)~|\frac12\Delta_k\hat\varphi_p(l,z)|}{K
 \sum_{(j,j')=(0,\pm1),(1,-1)}|\hat G_{\mu_p(z)}(l+jk)\,\hat G_{\mu_p
 (z)}(l+j'k)|}\le1+O(\lambda)\le2.
\end{align}
This completes the proof of Proposition~\ref{prop:bootstrapping}(ii). \QED

\subsection{Proof of Proposition~\ref{prop:bootstrapping}(i)}\label{ss:IR(i)}
Proposition~\ref{prop:bootstrapping}(i) is an immediate consequence of
Proposition~\ref{prop:diagbd}(i) and the following lemma:

\begin{lemma}\label{lem:WmTm-bd}
Let $d>2(\alpha\wedge2)$ and $L\gg1$, and fix $p<\pc$ and $m<m_p$.  Then,
$f(p,m)\le3$ implies that there are $(p,m)$-independent constants $C_{\sss
T},C_{\sss W}<\infty$ such that
\begin{align}\lbeq{WmTm-bd}
 T_p^{\sss(m)}\le C_{\sss T}\lambda,&&
 W_p^{\sss(m)}(k)\le C_{\sss W}\lambda\hat G_{\mu_p(m\vee1)}(k)^{-1}.
\end{align}
\end{lemma}

\Proof{Proof}
 Note that the Fourier transform of
$\varphi_p^{\sss(m)}(x,n)\equiv\varphi_p(x,n)m^n$ for $m<m_p$ is
\begin{align}
\hat\varphi_p^{\sss(m)}(k,e^{i\theta})=\sum_{(x,n)}\varphi_p^{\sss(m)}(x,n)
 e^{ik\cdot x}e^{i\theta n}=\sum_{(x,n)}\varphi_p(x,n)e^{ik\cdot x}(me^{i
 \theta})^n=\hat\varphi_p(k,me^{i\theta}).
\end{align}
By $f_1(p,m)\vee f_2(p,m)\le3$ and \refeq{G-IRbd}, $T_p^{\sss(m)}$ is bounded
as
\begin{align}\lbeq{triangle}
T_p^{\sss(m)}&\le p^2m\int_{[-\pi,\pi]^d}\frac{d^dk}{(2\pi)^d}\,\hat D
 (k)^2\int_{-\pi}^\pi\frac{d\theta}{2\pi}\,|\hat\varphi_p(k,e^{i
 \theta})|^2|\hat\varphi_p^{\sss(m)}(k,e^{-i\theta})|\nn\\
&~\le3^2\int_{[-\pi,\pi]^d}\frac{d^dk}{(2\pi)^d}\,\hat D(k)^2\int_{-\pi}^\pi
 \frac{d\theta}{2\pi}\bigg(\frac{3c}{\frac1{\hat\varphi_p(0,1)}+|\theta|+1-
 \hat D(k)}\bigg)^2\frac{3c}{\frac1{\hat\varphi_p(0,m)}+|\theta|+1-\hat D
 (k)}\nn\\
&~\le O(1)\int_{[-\pi,\pi]^d}\frac{d^dk}{(2\pi)^d}\,\frac{\hat D(k)^2}{(1-\hat
 D(k))^2}=O(1)\sum_{n=2}^\infty(n-1)\,D^{\star n}(o)\le O(\lambda),
\end{align}
where the last inequality is due to \refeq{Dprop} and $d>2(\alpha\wedge2)$.

To prove the bound on $W_p^{\sss(m)}(k)$, we first note that, by
$(q_p*\varphi_p)(y,t)\le(q_p*q_p*\varphi_p)(y,t)$ for $t\ge2$,
\begin{align}\lbeq{Wmbd-dec}
&\sum_{(y,t)}\big(1-\cos(k\cdot y)\big)(q_p*\varphi_p)(y,t)\cdot(q_p*
 \varphi_p)(y-x,t-n)\nn\\
&=p\sum_{y\in\Zd}\big(1-\cos(k\cdot y)\big)D(y)\cdot(q_p*\varphi_p)(y-x,
 1-n)\nn\\
&\quad+\sum_{(y,t):t\ge2}\big(1-\cos(k\cdot y)\big)(q_p*q_p*\varphi_p)(y,
 t)\cdot(q_p*\varphi_p)(y-x,t-n).
\end{align}
In the first sum on the right-hand side of \refeq{Wmbd-dec}, $1-n$ must be
larger than or equal to 1. If $1-n=1$, then, since
$(q_p*\varphi_p)(y-x,1)\equiv pD(y-x)\le p\|D\|_\infty\le Cp\lambda$ (see
\refeq{Dprop}), $f_1(p,m)\le3$ and $1-\hat D(k)\le2\hat
G_{\mu_p(m\vee1)}(k)^{-1}$ (see below \refeq{Delta1-a}), we obtain
\begin{align}
p\sum_y\big(1-\cos(k\cdot y)\big)D(y)\cdot(q_p*\varphi_p)(y-x,1)m\le3^2C
 \lambda\big(1-\hat D(k)\big)\le18C\lambda\hat G_{\mu_p(m\vee1)}(k)^{-1}.
\end{align}
If $1-n\ge2$, then we use
$(q_p*\varphi_p)(y,1-n)\le(q_p*q_p*\varphi_p)(y,1-n)$, $f_1(p,m)\vee
f_2(p,m)\le3$, \refeq{G-IRbd} and $1-\hat D(k)\le2\hat
G_{\mu_p(m\vee1)}(k)^{-1}$ to obtain that, for $m<1$,
\begin{align}
&p\sum_y\big(1-\cos(k\cdot y)\big)D(y)\cdot(mq_p*mq_p*\varphi_p^{\sss(m)})
 (y-x,1-n)\nn\\
&\le3\big(1-\hat D(k)\big)\int_{[-\pi,\pi]^d}\frac{d^dl}{(2\pi)^d}\,\hat D
 (l)^2\int_{-\pi}^\pi\frac{d\theta}{2\pi}\,\frac{3^3c}{\frac1{\hat\varphi_p
 (0,m)}+|\theta|+1-\hat D(l)}\nn\\
&\le O(1)\,\hat G_{\mu_p(m\vee1)}(k)^{-1}\int_{[-\pi,\pi]^d}\frac{d^dl}{(2
 \pi)^d}\,\frac{\hat D(l)^2}{1-\hat D(l)}\le O(\lambda)\,\hat G_{\mu_p(m
 \vee1)}(k)^{-1},
\end{align}
where the last inequality is due to \refeq{Dprop} and $d>\alpha\wedge2$.  The
other case of $m\ge1$ can be estimated in the same way.

To complete the proof of the bound on $W_p^{\sss(m)}(k)$, it remains to show
that the second sum on the right-hand side of \refeq{Wmbd-dec} is bounded by a
multiple of $\lambda\hat G_{\mu_p(m\vee1)}(k)^{-1}$.  Using
$1-\cos\sum_{j=1}^3X_j\le7\sum_{j=1}^3(1-\cos X_j)$ (cf.,
\cite[(4.50)]{Slad06}), we have
\begin{align}
\big(1-\cos(k\cdot y)\big)(q_p*q_p*\varphi_p)(y,t)&\le7p^2\sum_{u,v\in
 \Zd}\bigg(\big(1-\cos(k\cdot u)\big)D(u)D(v-u)\varphi_p(y-v,t-2)\nn\\
&\qquad+D(u)\Big(1-\cos\big(k\cdot(v-u)\big)\Big)D(v-u)\varphi_p(y-v,t-2)
 \nn\\
&\qquad+D(u)D(v-u)\Big(1-\cos\big(k\cdot(y-v)\big)\Big)\varphi_p(y-v,t-2)
 \bigg).
\end{align}
Recalling \refeq{DeltaFourier} and using $f_1(p,m)\le3$, we obtain that, for
$m<1$,
\begin{align}\lbeq{1-cos.phi-dec}
&\sum_{(y,t)}\big(1-\cos(k\cdot y)\big)(q_p*q_p*\varphi_p)(y,t)\cdot(mq_p*
 \varphi_p^{\sss(m)})(y-x,t-n)\nn\\
&\le7\cdot3^3\bigg(2\big(1-\hat D(k)\big)\int_{[-\pi,\pi]^d}\frac{d^dl}{(2
 \pi)^d}\,\hat D(l)^2\int_{-\pi}^\pi\frac{d\theta}{2\pi}\,|\hat\varphi_p(l,
 e^{i\theta})\,\hat\varphi_p^{\sss(m)}(l,e^{-i\theta})|\nn\\
&\qquad\qquad+\int_{[-\pi,\pi]^d}\frac{d^dl}{(2\pi)^d}\,|\hat D(l)|^3\int_{
 -\pi}^\pi\frac{d\theta}{2\pi}\,\big|\tfrac12\Delta_k\hat\varphi_p(l,e^{i
 \theta})\big||\hat\varphi_p^{\sss(m)}(l,e^{-i\theta})|\bigg).
\end{align}
Similarly to the above, by using $1-\hat D(k)\le2\hat G_{\mu_p(m\vee1)}
(k)^{-1}$, $f_2(p,m)\le3$ and \refeq{G-IRbd}, the first term on the right-hand
side of \refeq{1-cos.phi-dec} is bounded by a multiple of $\lambda\hat
G_{\mu_p(m\vee1)}(k)^{-1}$ when $d>\alpha\wedge2$.  For the second term on the
right-hand side of \refeq{1-cos.phi-dec}, we use $f_2(p,m)\vee f_3(p,m)\le3$ to
obtain
\begin{align}\lbeq{Wmbd-last}
&\int_{[-\pi,\pi]^d}\frac{d^dl}{(2\pi)^d}\,|\hat D(l)|^3\int_{-\pi}^\pi\frac{
 d\theta}{2\pi}\,\big|\tfrac12\Delta_k\hat\varphi_p(l,e^{i\theta})\big||\hat
 \varphi_p^{\sss(m)}(l,e^{-i\theta})|\nn\\
&\le3^2K\hat G_{\mu_p(m\vee1)}(k)^{-1}\sum_{(j,j')=(0,\pm1),(1,-1)}\int_{[-
 \pi,\pi]^d}\frac{d^dl}{(2\pi)^d}\,\hat D(l)^2\nn\\
&\hspace{10pc}\times\int_{-\pi}^\pi\frac{d\theta}{2\pi}\,|\hat G_{\mu_p(e^{i
 \theta})}(l+jk)||\hat G_{\mu_p(e^{i\theta})}(l+j'k)||\hat G_{\mu_p(me^{-i
 \theta})}(l)|.
\end{align}
By using \refeq{G-IRbd} as in \refeq{triangle}, the summand is bounded by a
multiple of $\lambda$ for any $k$ (the worst case is when $k=0$) as long as
$d>2(\alpha\wedge2)$.  This completes the proof of the bound on
$W_p^{\sss(m)}(k)$ and of Lemma~\ref{lem:WmTm-bd}. \QED

\section{Proof of Proposition~\ref{prop:piderbd}}\label{s:pc}
In this section, we prove Proposition~\ref{prop:piderbd} that was used in
Section~\ref{ss:IRbound} to prove Theorem~\ref{thm:pc}.  First we note that, by
\refeq{pider-exp},
\begin{align}
|\partial_p\hat\pi_p(0,1)|\le\frac1p\sum_{N=1}^\infty\sum_{(x,n)}
 \Pi_p^{\sss(N)}(x,n),
\end{align}
where $\Pi_p^{\sss(N)}(x,n)$ obeys the diagrammatic bound \refeq{PiN-diagbd},
with $T_p^{\sss(1)}\le C_{\sss T}\lambda$ as in \refeq{WmTm-bd}. Therefore, to
complete the proof of Proposition~\ref{prop:piderbd}, it suffices to prove the
following lemma:

\begin{lemma}\label{lem:tildeTHbd}
Let $d>2(\alpha\wedge2)$ and $L\gg1$.  Then, there are $C_{\sss\tilde
T},C_{\sss H}<\infty$ such that, for $p\in(1,\pc)$,
\begin{align}\lbeq{tildeTHbd}
\tilde T_p\le C_{\sss\tilde T}\lambda,&& H_p\le C_{\sss H}\lambda^2.
\end{align}
\end{lemma}

\Proof{Proof} The bound on $\tilde T_p$ can be proved in the same way as in
\refeq{triangle}.  Taking the Fourier transform and using
Theorem~\ref{thm:IRbound} and $\pc=1+O(\lambda)\le2$, we can bound $H_p$ as
\begin{align}\lbeq{H-1stbd}
H_p&\le p^5\int_{[-\pi,\pi]^{2d}}\frac{d^dk_1}{(2\pi)^d}\frac{d^dk_2}{(2
 \pi)^d}\,\hat D(k_1)^2\hat D(k_2)^2\big|\hat D(k_1-k_2)\big|\int_{-\pi}^\pi
 \frac{d\theta_1}{2\pi}\,\big|\hat\varphi_p(k_1,e^{i\theta_1})\big|^2\nn\\
&\hspace{7pc}\times\int_{-\pi}^\pi\frac{d\theta_2}{2\pi}\,\big|\hat\varphi_p
 (k_2,e^{i\theta_2})\big|^2\big|\hat\varphi_p(k_1-k_2,e^{i(\theta_1-\theta_2
 )})\big|^2\nn\\
&\le2^5\int_{[-\pi,\pi]^{2d}}\frac{d^dk_1}{(2\pi)^d}\frac{d^dk_2}{(2\pi)^d}
 \,\hat D(k_1)^2\hat D(k_2)^2\int_{-\pi}^\pi\frac{d\theta_1}{2\pi}\bigg(
 \frac{C}{|\theta_1|+1-\hat D(k_1)}\bigg)^2\nn\\
&\hspace{7pc}\times\int_{-\pi}^\pi\frac{d\theta_2}{2\pi}\bigg(\frac{C}{|
 \theta_2|+1-\hat D(k_2)}\bigg)^2\bigg(\frac{C}{|\theta_1-\theta_2|+1-\hat
 D(k_1-k_2)}\bigg)^2\nn\\
&\le\int_{[-\pi,\pi]^{2d}}\frac{d^dk_1}{(2\pi)^d}\frac{d^dk_2}{(2\pi)^d}
 \,\frac{\hat D(k_1)^2\hat D(k_2)^2}{(1-\hat D(k_2))^2}\int_{-\pi}^\pi\frac
 {d\theta_1}{2\pi}\,\frac{O(1)}{(|\theta_1|+1-\hat D(k_1))^2(|
 \theta_1|+1-\hat D(k_1-k_2))}\nn\\
&\le O(1)\int_{[-\pi,\pi]^d}\frac{d^dk_2}{(2\pi)^d}\,\frac{\hat D(k_2)^2}{(
 1-\hat D(k_2))^2}\int_{[-\pi,\pi]^d}\frac{d^dk_1}{(2\pi)^d}\frac{\hat D
 (k_1)^2}{(1-\hat D(k_1)\vee\hat D(k_1-k_2))^2}.
\end{align}
To prove the bound on $H_p$ in \refeq{tildeTHbd}, it suffices to show that the
last integral with respect to $k_1$ is $O(\lambda)$ for every $k_2$.  Since
this is trivial if $\hat D(k_1)\ge\hat D(k_1-k_2)$ (then the integrals in
\refeq{H-1stbd} are decoupled, each of them is $O(\lambda)$), it is sufficient
to prove that
\begin{align}
\int_{[-\pi,\pi]^d}\frac{d^dk_1}{(2\pi)^d}\,\frac{\hat D(k_1)^2}{(1-\hat
 D(k_1-k_2))^2}=O(\lambda).
\end{align}
However, by \refeq{Dprop}, the integral over $\|k_1-k_2\|_\infty>(\ell L)^{-1}$
is bounded as
\begin{align}
\int_{\|k_1-k_2\|_\infty>(\ell L)^{-1}}\frac{d^dk_1}{(2\pi)^d}\,\frac{\hat
 D(k_1)^2}{(1-\hat D(k_1-k_2))^2}\le\frac1{\Delta^2}\int_{[-\pi,\pi]^d}
 \frac{d^dk_1}{(2\pi)^d}\,\hat D(k_1)^2\le\frac{\|D\|_\infty}{\Delta^2}
 =O(\lambda).
\end{align}
Moreover, by \refeq{Dprop-k0}, the integral over $\|k_1-k_2\|_\infty\le(\ell
L)^{-1}$ is bounded as, for $\alpha\ne2$,
\begin{align}\lbeq{pc-eval}
\int_{\|k_1-k_2\|_\infty\le(\ell L)^{-1}}\frac{d^dk_1}{(2\pi)^d}\,\frac
 {\hat D(k_1)^2}{(1-\hat D(k_1-k_2))^2}\le O(L^{-2(\alpha\wedge2)})\int_0
 ^{(\ell L)^{-1}}dr\,r^{d-1-2(\alpha\wedge2)}=O(\lambda).
\end{align}
The case for $\alpha=2$ can be estimated similarly, since the $\log$ divergence
as $|k|\to0$ in \refeq{Dprop-k0} is unimportant in \refeq{pc-eval} as long as
$d>4$. This completes the proof of Lemma~\ref{lem:tildeTHbd}. \QED

\section{Proof of Proposition~\ref{prp:tauber-suff}}\label{s:tauber}
In this section, we prove Proposition~\ref{prp:tauber-suff} that was used in
Section~\ref{ss:limiting} to show \refeq{limit-main}, the key ingredient for
the proof of Theorem~\ref{thm:limiting}.

First we derive an expression for $\partial_\zeta\hat\Phi_p(k,m_p\zeta)$. Since
$\hat A_p(k,z)=\hat A_p^{\sss(1)}(k)+\hat A_p^{\sss(2)}(k,z)$, where $\hat
A_p^{\sss(1)}(k)$ is independent of $z$, we have
\begin{align}\lbeq{psider-expr}
\partial_\zeta\hat\Phi_p(k,m_p\zeta)&\equiv\partial_\zeta\frac{-(1-\zeta)\hat
 A_p^{\sss(2)}(k,m_p\zeta)}{\big((1-\zeta)\hat A_p(k,m_p\zeta)+\hat B_p(k)
 \big)\big((1-\zeta)\hat A_p^{\sss(1)}(k)+\hat B_p(k)\big)}\nn\\
&=\frac{\hat A_p^{\sss(2)}(k,m_p\zeta)-(1-\zeta)\partial_\zeta\hat A_p^{\sss
 (2)}(k,m_p\zeta)}{\big((1-\zeta)\hat A_p(k,m_p\zeta)+\hat B_p(k)\big)\big(
 (1-\zeta)\hat A_p^{\sss(1)}(k)+\hat B_p(k)\big)}\nn\\
&\quad+\frac{-(1-\zeta)\hat A_p^{\sss(2)}(k,m_p\zeta)}{(1-\zeta)\hat A_p^{
 \sss(1)}(k)+\hat B_p(k)}~\frac{\hat A_p^{\sss(1)}(k)+\hat A_p^{\sss(2)}(k,
 m_p\zeta)-(1-\zeta)\partial_\zeta\hat A_p^{\sss(2)}(k,m_p\zeta)}{\big((1-
 \zeta)\hat A_p(k,m_p\zeta)+\hat B_p(k)\big)^2}\nn\\
&\quad+\frac{-(1-\zeta)\hat A_p^{\sss(2)}(k,m_p\zeta)}{(1-\zeta)\hat A_p(k,
 m_p\zeta)+\hat B_p(k)}~\frac{\hat A_p^{\sss(1)}(k)}{\big((1-\zeta)\hat A_p
 ^{\sss(1)}(k)+\hat B_p(k)\big)^2}.
\end{align}
Recall that
\begin{align}\lbeq{A2-recall}
\hat A_p^{\sss(2)}(k,m_p\zeta)=\frac{\partial_\zeta\hat\pi_p(k,m_p)^{-1}}{p
 m_p}-\frac{\hat\pi_p(k,m_p)^{-1}-\hat\pi_p(k,m_p\zeta)^{-1}}{pm_p(1-\zeta)},
\end{align}
where $\partial_\zeta\hat\pi_p(k,m_p)^{-1}$ is an abbreviation of
$\partial_\zeta\hat\pi_p(k,m_p \zeta)^{-1}|_{\zeta=1}\equiv
m_p\partial_z\hat\pi_p(k,z)^{-1}|_{z=m_p}$, so that
\begin{align}
(1-\zeta)\,\partial_\zeta\hat A_p^{\sss(2)}(k,m_p\zeta)&=(1-\zeta)\,
 \partial_\zeta\bigg(-\frac{\hat\pi_p(k,m_p)^{-1}-\hat\pi_p(k,m_p\zeta)^{-1}}
 {pm_p(1-\zeta)}\bigg)\nn\\
&=\frac{\partial_\zeta\hat\pi_p(k,m_p\zeta)^{-1}}{pm_p}-\frac{\hat\pi_p(k,
 m_p)^{-1}-\hat\pi_p(k,m_p\zeta)^{-1}}{pm_p(1-\zeta)}\nn\\
&=\frac{\partial_\zeta\hat\pi_p(k,m_p\zeta)^{-1}-\partial_\zeta\hat\pi_p(k,
 m_p)^{-1}}{pm_p}+\hat A_p^{\sss(2)}(k,m_p\zeta).
\end{align}
Therefore, $\hat A_p^{\sss(2)}(k,m_p\zeta)-(1-\zeta)\partial_\zeta\hat
A_p^{\sss(2)}(k,m_p\zeta)$ in \refeq{psider-expr} can be replaced by $\hat
a_p^{\sss(2)}(k,m_p\zeta)$, which is
\begin{align}\lbeq{a2-def}
\hat a_p^{\sss(2)}(k,m_p\zeta)=\frac{\partial_\zeta\hat\pi_p(k,m_p)^{-1}
 -\partial_\zeta\hat\pi_p(k,m_p\zeta)^{-1}}{pm_p}.
\end{align}

Now, Proposition~\ref{prp:tauber-suff} is an immediate consequence of the
following lemma:

\begin{lemma}\label{lem:psiingr-bd}
Let $d>2(\alpha\wedge2)$,
$\epsilon\in(0,1\wedge\frac{d-2(\alpha\wedge2)}{\alpha\wedge2})$ and $L\gg1$.
Then, the following hold uniformly in $p\in(0,\pc]$, $k\in[-\pi,\pi]^d$ and
$\zeta\in\mC$ with $|\zeta|<1$:
\begin{enumerate}[(i)]
\item
There is a positive constant $c$ such that
\begin{align}\lbeq{psiingr-bd1}
\left.\begin{array}{l}
 |(1-\zeta)\hat A_p(k,m_p\zeta)+\hat B_p(k)|\\
 |(1-\zeta)\hat A_p^{\sss(1)}(k)+\hat B_p(k)|
\end{array}\right\}\ge c|1-\zeta|.
\end{align}
\item
There is a finite constant $c_\epsilon$ such that
\begin{align}\lbeq{psiingr-bd2}
\left.\begin{array}{l}
 |A_p^{\sss(2)}(k,m_p\zeta)|\\
 |\hat a_p^{\sss(2)}(k,m_p\zeta)|
\end{array}\right\}\le c_\epsilon|1-\zeta|^\epsilon.
\end{align}
\end{enumerate}
\end{lemma}

In the following proof, the constant in the $O(\,\cdot\,)$ term is independent
of $p$, $k$ and $\zeta$.

\Proof{Proof of Lemma~\ref{lem:psiingr-bd}(i)}
 Since both bounds can be proved in the same way, we only prove the bound on
$|(1-\zeta)\hat A_p^{\sss(1)}(k)+\hat B_p(k)|$.

We consider the following four cases: (a)~$\Re\zeta\le0$; (b)~$\Re\zeta\ge0$
with $\Re(1-\zeta)\ge\Im(1-\zeta)$; (c)~$\Re\zeta\ge0$ with
$\Re(1-\zeta)\le\Im(1-\zeta)$ and $\hat D(k)\ge1-\Delta$; (d)~$\Re\zeta\ge0$
and $\hat D(k)\le1-\Delta$.  Note that these four cases exhaust all
$\zeta\in\mC$ with $|\zeta|\le1$.  For the moment, we abbreviate $\hat
A_p^{\sss(1)}(k)$ to $A$ and $\hat B_p(k)$ to $B$.
\begin{enumerate}[(a)]
\item
Since $A,B\in\mR$ and $|w|\ge|\Re w|$ for any $w\in\mC$,
\begin{align}
|(1-\zeta)A+B|=|A+B-A\zeta|\ge|A+B+A\Re(-\zeta)|.
\end{align}
Since $\hat D(k)>-1+\Delta$ holds for all $k\in[-\pi,\pi]^d$ (cf.,
\refeq{Dprop}), we have $A=\hat D(k)+O(\lambda)\ge-1+\Delta-O(\lambda)$. Since
$A+B=1+O(\lambda)$ and $\Re(-\zeta)\ge0$, we obtain
\begin{align}\lbeq{Rez<0}
|A+B+A\Re(-\zeta)|=A+B+A\Re(-\zeta)\ge\Delta-O(\lambda),
\end{align}
uniformly in the concerned $\zeta$.
\item[(d)]
Using $\Re\zeta\ge0$ and $\hat D(k)\le1-\Delta$, we can prove \refeq{Rez<0}
similarly.
\item
Since $A+B=1+O(\lambda)$, $B\ge0$, $\Re\zeta\ge0$, and
$\Re(1-\zeta)\ge\frac1{\sqrt2}|1-\zeta|$ for $\zeta$ in case~(ii), we obtain
\begin{align}
|(1-\zeta)A+B|=|(1-\zeta)(A+B)+B\zeta|&\ge|(A+B)\,\Re(1-\zeta)+B\Re\zeta|\nn\\
&\ge(A+B)\,\Re(1-\zeta)\ge\frac{1-O(\lambda)}{\sqrt2}\,|1-\zeta|.
\end{align}
\item
Since $A=\hat D(k)+O(\lambda)\ge1-\Delta-O(\lambda)$ and
$|\Im(1-\zeta)|\ge\frac1{\sqrt2}|1-\zeta|$ for $\zeta$ in case~(iii), by using
the imaginary part (i.e., $|w|\ge|\Im w|$ for $w\in\mC$) we obtain
\begin{align}
|(1-\zeta)A+B|\ge|A\Im(1-\zeta)|=A|\Im(1-\zeta)|\ge\frac{1-\Delta-O(\lambda)}
 {\sqrt2}\,|1-\zeta|.
\end{align}
\end{enumerate}
This completes the proof of Lemma~\ref{lem:psiingr-bd}(i).
\QED

\Proof{Proof of Lemma~\ref{lem:psiingr-bd}(ii)}
 First, by adding and subtracting, we can rewrite
$\hat A_p^{\sss(2)}(k,m_p\zeta)$ in \refeq{A2-recall} as
\begin{align}\lbeq{A2-rewr}
\hat A_p^{\sss(2)}(k,m_p\zeta)&=-\frac{\partial_\zeta\hat\pi_p(k,m_p)}{pm_p
 \hat\pi_p(k,m_p)^2}+\frac{\hat\pi_p(k,m_p)-\hat\pi_p(k,m_p\zeta)}{pm_p(1-
 \zeta)\hat\pi_p(k,m_p)\hat\pi_p(k,m_p\zeta)}\nn\\
&=-\frac{\partial_\zeta\hat\pi_p(k,m_p)-\frac{\hat\pi_p(k,m_p)-\hat\pi_p(k,
 m_p\zeta)}{1-\zeta}}{pm_p\hat\pi_p(k,m_p)^2}+\frac{\big(\hat\pi_p(k,m_p)-
 \hat\pi_p(k,m_p\zeta)\big)^2}{pm_p(1-\zeta)\hat\pi_p(k,m_p)^2\hat\pi_p(k,
 m_p\zeta)},
\end{align}
and $\hat a_p^{\sss(2)}(k,m_p\zeta)$ in \refeq{a2-def} as
\begin{align}\lbeq{a2-rewr}
\hat a_p^{\sss(2)}(k,m_p\zeta)&=-\frac{\partial_\zeta\hat\pi_p(k,m_p)}{pm_p
 \hat\pi_p(k,m_p)^2}+\frac{\partial_\zeta\hat\pi_p(k,m_p\zeta)}{pm_p\hat
 \pi_p(k,m_p\zeta)^2}\nn\\
&=-\frac{\partial_\zeta\hat\pi_p(k,m_p)-\partial_\zeta\hat\pi_p(k,m_p\zeta)}
 {pm_p\hat\pi_p(k,m_p)^2}+\frac{\big(\hat\pi_p(k,m_p)^2-\hat\pi_p(k,m_p
 \zeta)^2\big)\,\partial_\zeta\hat\pi_p(k,m_p\zeta)}{pm_p\hat\pi_p(k,m_p)^2
 \hat\pi_p(k,m_p\zeta)^2}.
\end{align}
Since $pm_p\hat\pi_p(k,m_p)^2=1+O(\lambda)$,
$\hat\pi_p(k,m_p\zeta)^2=1+O(\lambda)$,
$|\partial_\zeta\hat\pi_p(k,m_p\zeta)|=O(\lambda)$ and
$|\hat\pi_p(k,m_p)-\hat\pi_p(k,m_p\zeta)|=O(\lambda)|1-\zeta|$, the second
terms in \refeq{A2-rewr}--\refeq{a2-rewr} are $O(|1-\zeta|)$.  To prove
\refeq{psiingr-bd2}, it thus suffices to show that the numerator of the first
term in \refeq{A2-rewr} and that in \refeq{a2-rewr} are both bounded by
$O_\epsilon(1)|1-\zeta|^\epsilon$, where the constant in the $O_\epsilon(1)$
term may depend on $\epsilon$.  Since both can be proved similarly, we only
prove that $|\partial_\zeta\hat\pi_p(k,m_p)-\partial_\zeta\hat
\pi_p(k,m_p\zeta)|\le O_\epsilon(1)|1-\zeta|^\epsilon$.

Note that
\begin{align}\lbeq{fracder-time0}
|\partial_\zeta\hat\pi_p(k,m_p)-\partial_\zeta\hat\pi_p(k,m_p\zeta)|&=\bigg|
 \sum_{(x,n):n\ge2}n(1-\zeta^{n-1})\pi_p(x,n)e^{ik\cdot x}m_p^n\bigg|\nn\\
&\le\sum_{(x,n):n\ge2}n|1-\zeta^{n-1}|\,|\pi_p(x,n)|m_p^n.
\end{align}
For $n\ge2$, we have
\begin{align}
|1-\zeta^{n-1}|=\bigg|(1-\zeta^{n-1})^{1-\epsilon}\bigg(\frac{1-\zeta^{n-1}}
 {1-\zeta}\bigg)^\epsilon(1-\zeta)^\epsilon\bigg|\le2^{1-\epsilon}\bigg|
 \sum_{l=0}^{n-2}\zeta^l\bigg|^\epsilon|1-\zeta|^\epsilon\le2|1-\zeta|
 ^\epsilon n^\epsilon.
\end{align}
Moreover, for $\epsilon\in(0,1)$, we have (cf., \cite[Section~6.3]{ms93})
\begin{align}
n^{1+\epsilon}=\frac{n^2}{(1-\epsilon)\,\Gamma(1-\epsilon)}\int_0^\infty
 e^{-n\rho^{1/(1-\epsilon)}}d\rho.
\end{align}
Applying these to \refeq{fracder-time0} and using the diagrammatic bound
\refeq{piNn-diagbd} for $r=2$ and $T_p^{\sss(\tilde m_\rho)}\!\le C_{\sss
T}\lambda$ with $\tilde m_\rho=m_pe^{-\rho^{1/(1-\epsilon)}}$, we have
\begin{align}\lbeq{fracder-time}
|\partial_\zeta\hat\pi_p(k,m_p)-\partial_\zeta\hat\pi_p(k,m_p\zeta)|&\le\frac{
 2|1-\zeta|^\epsilon}{(1-\epsilon)\,\Gamma(1-\epsilon)}\int_0^\infty d\rho~
 \sum_{(x,n)}n^2|\pi_p(x,n)|\,\tilde m_\rho^n\\
&\le\frac{2(1+2C_{\sss T}\lambda)|1-\zeta|^\epsilon}{(1-\epsilon)\,\Gamma(1-
 \epsilon)}\sum_{N=0}^\infty(N+1)^2(2C_{\sss T}\lambda)^{(N-1)\vee0}\int_0
 ^\infty d\rho~S_p^{\sss(\tilde m_\rho)},\nn
\end{align}
where, by \refeq{IRbound} and $p\le pm_p=1+O(\lambda)$,
\begin{align}
\int_0^\infty d\rho~S_p^{\sss(\tilde m_\rho)}&\le p^2m_p\int\frac{d^dk}
 {(2\pi)^d}\,\hat D(k)^2\int\frac{d\theta}{2\pi}\bigg(\frac{C}{p(m_p-1)
 +|\theta|+1-\hat D(k)}\bigg)^3\nn\\
&\hspace{10pc}\times\int_0^\infty d\rho\,\frac{Ce^{-\rho^{1/(1-\epsilon)
 }}}{pm_p(1-e^{-\rho^{1/(1-\epsilon)}})+|\theta|+1-\hat D(k)}\nn\\
&\le\int\frac{d^dk}{(2\pi)^d}\,\frac{O(1)}{(1-\hat D(k))^2}\int_0^\infty
 d\rho\,\frac{e^{-\rho^{1/(1-\epsilon)}}}{1-e^{-\rho^{1/(1-\epsilon)}}+1
 -\hat D(k)}.
\end{align}
However, since
\begin{align}
\int_0^\infty d\rho\,\frac{e^{-\rho^{1/(1-\epsilon)}}}{1-e^{-\rho^{1/(1-
 \epsilon)}}+1-\hat D(k)}&=\int_0^\infty ds\,\frac{(1-\epsilon)s^{-\epsilon}
 e^{-s}}{1-e^{-s}+1-\hat D(k)}\qquad(\because\rho=s^{1-\epsilon})\nn\\
&\le\frac{1-\epsilon}{1-e^{-1}}\bigg(\int_1^\infty ds\;e^{-s}+\int_0^1ds\,
 \frac{s^{-\epsilon}}{s+1-\hat D(k)}\bigg)\nn\\
&\le\frac{1-\epsilon}{1-e^{-1}}\bigg(1+\int_0^{1-\hat D(k)}ds\,\frac{s^{-
 \epsilon}}{1-\hat D(k)}+\int_{{1-\hat D(k)}}^1ds\;s^{-1-\epsilon}\bigg)\nn\\
&\le\frac{1-\epsilon}{1-e^{-1}}\bigg(1+\frac{(1-\hat D(k))^{-\epsilon}}{1
 -\epsilon}+\frac{(1-\hat D(k))^{-\epsilon}}\epsilon\bigg),
\end{align}
we obtain that
\begin{align}\lbeq{square-bd}
\int_0^\infty d\rho~S_p^{\sss(\tilde m_\rho)}\le\int\frac{d^dk}{(2\pi)^d}\,
 \frac{O_\epsilon(1)}{(1-\hat D(k))^{2+\epsilon}}<\infty,
\end{align}
as long as $d>(2+\epsilon)(\alpha\wedge2)$, due to \refeq{Dprop-k0}.  By
\refeq{fracder-time}, this completes the proof of
$|\partial_\zeta\hat\pi_p(k,m_p)-\partial_\zeta\hat\pi_p(k,m_p\zeta)|\le
O_\epsilon(1)|1-\zeta|^\epsilon$ and of Lemma~\ref{lem:psiingr-bd}(ii). \QED

\appendix
\section{Proof of Proposition~\ref{prop:Dnsup-bd}}\label{a:Dprop}
In this section, we prove the bounds on $D$ summarized in
Proposition~\ref{prop:Dnsup-bd}.  Since the bounds on $1-\hat D(k)$ in
\refeq{Dprop} are equivalent to \cite[(1.20)--(1.21)]{hs02} whose proofs are
independent of the range of $D$ (see \cite[Appendix~A]{hs02}), it thus remains
to prove the bound on $\|D^{\star n}\|_\infty$ in \refeq{Dprop} and the bounds
on $1-\hat D(k)$ for $\|k\|_\infty\le(\ell L)^{-1}$ in \refeq{Dprop-k0}.

First we prove the bound on $\|D^{\star n}\|_\infty$ assuming \refeq{Dprop-k0}.
By definition, it is trivial when $n=1$.  For $n\ge2$, we let
\begin{align}
R=\{k\in[-\pi,\pi]^d:|k|\le(\ell L)^{-1},~\hat D(k)\ge0\},
\end{align}
so that $|\hat D(k)|=1-(1-\hat D(k))\le e^{-(1-\hat D(k))}$ for $k\in R$, and
that $0\le|\hat D(k)|<1-\Delta$ for $k\notin R$, due to the bound on $1-\hat
D(k)$ in \refeq{Dprop}.  Therefore, for any $x\in\Zd$,
\begin{align}
D^{\star n}(x)\le\int_{[-\pi,\pi]^d}\frac{d^dk}{(2\pi)^d}\;|\hat D(k)|^n
 \le\int_R\frac{d^dk}{(2\pi)^d}\;e^{-n(1-\hat D(k))}+(1-\Delta)^{n-2}
 \int_{R^\text{c}}\frac{d^dk}{(2\pi)^d}\;\hat D(k)^2,
\end{align}
where the integral over $k\in R^\text{c}\equiv[-\pi,\pi]^d\setminus R$ is
bounded by $\|D\|_\infty(1-\Delta)^{n-2}\le O(\lambda)\,n^{-d/(\alpha\wedge
2)}$.  For the integral over $k\in R$, we use the bounds on $1-\hat D(k)$ in
\refeq{Dprop-k0}.  If $\alpha\ne2$, then
\begin{align}\lbeq{Dnsup-ne2}
\int_R\frac{d^dk}{(2\pi)^d}\;e^{-n(1-\hat D(k))}\le c'\lambda\int_0^\infty
 \frac{dr}r\;r^de^{-cnr^{\alpha\wedge2}}=\frac{c'\Gamma(\frac{d}{\alpha
 \wedge2})\lambda}{(\alpha\wedge2)(cn)^{d/(\alpha\wedge2)}},
\end{align}
for some $c,c'\in(0,\infty)$, where $r=\ell L|k|$.  If $\alpha=2$, then
\begin{align}\lbeq{Dnsup-e2}
\int_R\frac{d^dk}{(2\pi)^d}\;e^{-n(1-\hat D(k))}\le c'\lambda\int_0^1
 \frac{dr}r\;r^de^{-cnr^2\log\frac\pi{2r}}\le c'\lambda\int_0^\infty
 \frac{dr}r\;r^de^{-c''nr^2}=\frac{c'\Gamma(\frac{d}2)\lambda}{2(c''
 n)^{d/2}},
\end{align}
where $c''=c\log\frac\pi2>0$.  This completes the proof of the bound on
$\|D^{\star n}\|_\infty$ in \refeq{Dprop}.

Next we prove the bounds on $1-\hat D(k)$ for $|k|\le(\ell L)^{-1}$ with
$L\gg1$. Since $\|k\|_\infty\le|k|$, this is sufficient for the proof of
\refeq{Dprop-k0}. First we note that, by the Riemann sum approximation,
\begin{align}
\frac1{L^d}\sum_{x\in\Zd}h(x/L)=\int_{\Rd}d^dx\;h(x)+o(1)=1+o(1),
\end{align}
where $o(1)\to0$ as $L\to\infty$.  Therefore,
\begin{align}
1-\hat D(k)=\big(1+o(1)\big)(I_1+I_2+I_3),
\end{align}
where
\begin{align}
I_1&=L^{-d}\sum_{x\in\Zd:|x|<\ell L}h(x/L)\big(1-\cos(k\cdot x)\big),\\
I_2&=L^{-d}\sum_{\substack{x\in\Zd:\\ \ell L\le|x|<\frac\pi{2|k|}}}h(x/L)
 \big(1-\cos(k\cdot x)\big),\\
I_3&=L^{-d}\sum_{x\in\Zd:|x|\ge\frac\pi{2|k|}}h(x/L)\big(1-\cos(k\cdot x)
 \big).
\end{align}
However, by \refeq{majorass} and using $1-\cos(k\cdot x)\asymp|k|^2|x|^2$ if
$|x|\le\frac\pi{2|k|}$ and $1-\cos(k\cdot x)\le2$ otherwise, we obtain
\begin{align}
I_1&\le O(L^{-d}|k|^2)\sum_{x\in\Zd:|x|<\ell L}|x|^2=O((L|k|)^2),\\
I_2&\asymp O(L^\alpha|k|^2)\sum_{\substack{x\in\Zd:\\ \ell L\le|x|<\frac\pi
 {2|k|}}}|x|^{-d-\alpha+2}=
 \begin{cases}
 O((L|k|)^{\alpha\wedge2})&(\alpha\ne2),\\
 O\big((L|k|)^2\log\tfrac\pi{2\ell L|k|}\big)&(\alpha=2),
 \end{cases}\\
I_3&\le O(L^\alpha)\sum_{x\in\Zd:|x|\ge\frac\pi{2|k|}}|x|^{-d-\alpha}=O((L|k|)
 ^\alpha).
\end{align}
This completes the proof of \refeq{Dprop-k0}. \QED

\section*{Acknowledgements}
This work was supported in part by the Institute of Mathematics at Academia
Sinica in Taiwan, and in part by EURANDOM and the Department of Mathematics and
Computer Science at TU/e in the Netherlands.  The work of LCC was also
supported in part by National Science Council, and the work of AS was also
supported in part by the Netherlands Organisation for Scientific Research and
in part by the London Mathematical Sociaety.  LCC would like to thank EURANDOM
for its hospitality during the visit in the period of August~3--27, 2005.  AS
would like to thank the Institute of Mathematics at Academia Sinica for a
comfortable and stimulating environment during the visits in the period of
November~27--December~17, 2005, and in the period of December~18,
2006--January~6, 2007.  We would like to thank Wei-Shih Yang and Narn-Rueih
Shieh for valuable conversations, and Markus Heydenreich, Mark Holmes and Remco
van der Hofstad for useful comments on the previous version of the manuscript.

\end{document}